\documentclass[final,1p,times]{elsarticle}

\usepackage{xcolor}
\usepackage[notref,notcite]{showkeys}
\usepackage{lineno,hyperref}
\usepackage{amsmath,amsfonts}
\usepackage{verbatim}
\usepackage{latexsym}
\usepackage{graphicx}
\usepackage{epstopdf}
\usepackage{subfig}
\usepackage{caption}
\usepackage{dsfont}
\usepackage{indentfirst}
\usepackage{booktabs}
\newtheorem{theorem}{Theorem}[section]
\newtheorem{lemma}[theorem]{Lemma}

\newtheorem{assumption}[theorem]{Assumption}
\newtheorem{example}[theorem]{Example}

\newtheorem{remark}[theorem]{Remark}

\newcommand{\lf}{\lfloor}
\newcommand{\rf}{\rfloor}
\def\lan{\langle} \def\ran{\rangle}
\def\d{\delta}
\newcommand{\D}{\Delta}

\def\eps{\varepsilon}
\def\e{\epsilon}
\def\th{\theta}
\def\ve{\vee}
\def\we{\wedge}
\def\a{\alpha} \def\g{\gamma}
\def\b{\beta}

\def\k{\kappa} \def\l{\lambda}

\def\K{\times}

\def\x{\xi}
\def\f{\forall}

\newcommand{\Ito}{It\^{o} }
\newcommand{\Holder}{H\"{o}lder }
\newcommand{\BDG}{Burkholder-Davis-Gundy }
\newcommand{\pr}{\textbf{Proof.} }
\newcommand{\Cheb}{Chebyshev }

\newcommand{\E}{\mathbb{E}}
\newcommand{\F}{\mathcal{F}}
\newcommand{\PP}{\mathbb{P}}
\newcommand{\II}{\mathbb{I}}

\newcommand{\R}{\mathbb{R}}

\newcommand{\fd}{f_{\Delta}}
\newcommand{\gd}{g_{\Delta}}
\newcommand{\tkk}{t_{k+1}}
\newcommand{\tk}{t_{k}}
\newcommand{\qu}{\quad}
\newcommand{\no}{\nonumber }
\newcommand{\intt}{\int_0^t}
\newcommand{\pd}{\pi_{\Delta}}

\newcommand{\yd}{Y_{\Delta} }
\newcommand{\TDfds}{f _{\Delta} (Z_1(s),Z_2(s)  )}
\newcommand{\TDgds}{g _{\Delta} (Z_1(s),Z_2(s)  )}

\newcommand{\ykk}{ y_{k+1}}
\newcommand{\ykd}{ y_{k- \delta_k}}
\newcommand{\yk}{y_{k}}

\makeatletter
\@addtoreset{equation}{section}
\makeatother

\modulolinenumbers[5]

\begin{document}

\begin{frontmatter}
%
\title{The truncated EM method for stochastic differential delay equations with variable delay}


\author[mainaddress]{Shounian Deng}
\author[address3]{Chen Fei }
\author[mainaddress]{Weiyin Fei \corref{correspondingauthor}}
\cortext[correspondingauthor]{Corresponding author}
\ead{wyfei@ahpu.edu.cn}
\author[thirdaryaddress]{Xuerong Mao }

\address[mainaddress]{Key Laboratory of Advanced Perception and Intelligent Control of High-end Equipment, Ministry of Education,  and  School of  Mathematics-Physics and Finance, Anhui Polytechnic University, Wuhu, Anhui 241000, China.}
\address[address3]{Business School, University of Shanghai for Science and Technology,  Shanghai  200093, China}
\address[thirdaryaddress]{Department of Mathematics and Statistics, University of Strathclyde, Glasgow G1 1XH, UK}

%
%

\begin{abstract}
This paper mainly investigates the strong convergence and stability   of the truncated Euler-Maruyama (EM)  method for stochastic differential delay equations with variable delay whose coefficients can be growing super-linearly.
By constructing  appropriate  truncated functions to control the super-linear growth of the original coefficients,
we present a type of the truncated EM method for such SDDEs with variable delay, which is proposed to be  approximated by
the value taken at the nearest grid points on the left of the delayed argument.
The strong convergence result (without order) of the method is established under the local Lipschitz plus generalized Khasminskii-type conditions and
the optimal  strong convergence  order $1/2$ can be obtained if the global monotonicity with $U$ function and polynomial growth conditions   are added to the assumptions. Moreover,
the partially  truncated EM method is proved to preserve the  mean-square and $H_\infty$ stabilities of the true solutions.
Compared with the known results on the truncated EM method for SDDEs,  a better order of strong convergence is obtained  under more relaxing conditions on the coefficients, and more refined technical estimates are developed  so as to overcome the challenges arising due to variable delay.
Lastly, some numerical examples are utilized to confirm the effectiveness of the theoretical results.
\end{abstract}


\end{frontmatter}

\linenumbers



\section{Introduction}

Stochastic differential delay equations (SDDEs)  play a significant part in many application fields, such as economy, finance, automatic control and population dynamics (see, e.g., \citep{Li_Mao2020AM,Mao2004delay_LV,Mao2004Delay_pop,Mao2005SDDE_pop}).
In general, SDDEs rarely have explicit solutions available, and therefore  one has to construct
appropriate numerical schemes  so as to approximate their solutions.
Moreover,
many important SDDE models often possess super-linear growth coefficients in practice, for example,
stochastic delay Lotka-Volterra model arising  in population dynamics of the form (\cite{Mao2004delay_LV})
\begin{align}\label{Delay_LV_model}
dX(t) = \textrm{diag} \Big (  X_1 (t), \cdots , X_d (t)  \Big )\Big [ ( b + AX(t - \tau)) dt
 + \sigma X(t) d B(t)\Big ],
\end{align}
where $B(t)$ is a Brownian motion and $\sigma = (\sigma _{ij})_{d \K d }$ is a matrix
representing the intensity of noise. If we apply the
explicit Euler-Maruyama (EM) scheme to the model \eqref{Delay_LV_model}, it can be shown  that this EM approximation may not converge  to the true solution in the strong mean-square sense  at a finite time point  (see, e.g.,\cite{Hutzenthaler2011PRS}).
As a result, this paper concerns the strong convergence analysis of numerical scheme  for such SDDEs with super-linearly growing drift and diffusion coefficients.

When the delay term vanishes, the underlying SDDEs reduce to the usual  stochastic differential equations (SDEs), numerical methods for which have been extensively investigated for the past decades under the global Lipschitz condition (e.g., \cite{Peter1992}, \cite{Higham_review}, \cite{Platen2010book}).
In 2012,   \citet{HJK2012AAP} firstly introduced an explicit method called  the tamed EM method to solve SDEs with super-linearly growing drift coefficient and  linearly growing diffusion coefficient. Since then, several explicit  schemes have been  proposed for   SDEs whose coefficients can be growing super-linearly under the local Lipschitz condition,
e.g., balanced EM method \cite{Zhang2013SIAM}\cite{ZhangMa2017APNUM}, truncated EM method \cite{Mao2015truncated,Mao2016rates}.
Recently, the attention of some researches was attracted  to the strong convergence of explicit numerical methods for super-linear SDEs with delay, i.e., SDDEs.
Motivated by the work of Mao \cite{Mao2015truncated,Mao2016rates},   
\citet{Guo2018NA} discussed the strong convergence of the truncated EM method for SDDEs under the local Lipschitz condition plus the generalized Khasminskii-type condition \eqref{Guo_cond}.
In a subsequent paper,
\citet{Hu_Gao2021FM} took the delay and jumps into consideration, they
extended  convergence results from \cite{Guo2018NA,Deng2019CAM}  to the case of SDDEs with Poisson jumps.
By using a different estimate  for the difference between  the original  and  the truncated coefficients,
 \citet{Fei2020CPAA} relaxed the restrictive condition   on the step size which
is required to extremely  small and thus improved the convergence results of \cite{Guo2018NA}.
%
Moreover,
\citet{SHGL2020NA} achieved a better convergence order  than \cite{Guo2018NA} and
\cite{Fei2020CPAA} by adopting the  truncation techniques  from \cite{LiMao2019IMA} for such  SDDEs.
The applications of the partially truncated  and modified truncated EM methods  for SDDEs can be found in \cite{Zhang_song2018CAM,Lan2019SDDE}.
Other explicit numerical  methods for super-linear SDDEs, say tamed EM, balanced EM, truncated Milstein, projected EM, are discussed in \cite{DKS2016SIAM,Deng2021CAM,Cao2021CAM,Song_Zhang2019CAM,Huang_Li2020CAM,Tan_Yuan2020SCI}.

In  the preceding discussion, the focus of the interest is on
the numerical method  of  SDEs with constant delay. In reality,  delay can behave as a function of time, namely, variable delay.
 \citet{Mao2003SDDE} were the first authors to  consider   the strong convergence of numerical  method for such SDDEs.
In contrast to the constant delay, the main difficulty  in the construction of the computational approach in case of SDDEs with variable delay   is how to numerically approximate  the values of the  solution at the delayed instants.
 \citeauthor{Mao2003SDDE} proposed to use the  approximate values at the nearest grid points on the left of the delayed arguments  to
estimate the variable delay,
 and they derived an upper bound for the difference of  this approximation.
Influenced by \cite{Mao2003SDDE}, several
explicit and implicit variants of the EM method have been developed
%
to study the convergence  of the numerical solutions to stochastic
equations with variable delay. For instance,  we refer to \cite{MM2011} for the convergence in probability of EM method for highly nonlinear neutral stochastic differential equations with variable delay, and to \cite{Deng2014HPC}
for the strong error analysis
of EM method for SDEs with variable and distributed delays.
%
Further, in \cite{Deng2019APNUM}, strong convergence rates are derived for the   split-step theta method applied to stochastic age-dependent population equations with Markovian switching and variable delay.

%
%
%
%
%
%
%
%

Based on the above discussion, the objective of this work is to study
 the strong convergence (including the stability)  of explicit numerical method of super-linear SDDEs with variable delay.
By constructing  appropriate  truncated functions to control the super-linear growth of the original coefficients,
we introduce a type of the  truncated EM method for such SDDEs with variable delay, which can be  approximated by  the value taken at the left endpoint of interval containing the delayed argument.
We then show that the method is convergent in the sense of mean-square according to the properties \eqref{PP2} and \eqref{PP3} of the truncated  functions:  preservation of  Khasminskii-type condition and  linear growth  for a fixed step size. In addition,
we discuss the mean-square and  $H_\infty$  stabilities   of the method.
\par

The main contribution of this paper is to improve the strong convergence results of the truncated EM method for SDDEs.
 We obtain a better  strong convergence order  than the existing results such as
\citet[Theorem 3.12]{GuoZhan2018IJCM},
\citet[Theorem 5.3]{Guo2018NA},
and
 \citet[Theorem 3.6]{Fei2020CPAA},
under the more relaxing conditions:  a generalized  Khasminskii-type   and
   global monotonicity conditions  with $U$ function, see Remarks \ref{rem1} and \ref{rem2}.
%
It should be pointed out that  the appearance  of $U$ function in global monotonicity
condition
will make the choice of the  coefficients for SDDEs more flexible, see e.g., \citet[Example 6.2]{Fei2020CPAA}.
Our technical estimates are more refined than those of \citet{Guo2018NA}, \citet{Fei2020CPAA} and \citet{SHGL2020NA}
  in that we develop new techniques to overcome the challenges arising due to variable delay.
%
%
Moreover, it is proved that the partially truncated EM method has the properties of
  the  mean-square and   $H_\infty$ stabilities.
\par
The remainder of the paper is organised as follows. The second section introduces some basic notions and assumptions. The next section describes the truncated EM scheme for SDDEs with variable delay.  Strong convergence results are established in the forth section. In the fifth section, we prove some
stability theorems. Numerical simulations are provided in the sixed section. In the final section, we close the paper by our conclusion.

\section{Preliminaries}\label{secmathpre}

Throughout this paper, let $(\Omega , {\mathcal F}, \mathbb{P})$ be a complete probability
space with a filtration  $\{{\cal F}_t\}_{t\ge 0}$ satisfying the usual conditions (i.e., it is increasing and right continuous while $\cal{F}_\textrm{0} $ contains all $\mathbb{P}$\textrm{-}null sets). Let $\tau > 0$ be a constant and denote by  $C([- \tau ,0 ]; \R ^d )$ the space of all continuous functions from $[ - \tau, 0 ]$ to $\R^d$ with
the norm $ \| \phi \| =\sup _ { -\tau \le  \th \le 0} | \phi (\th)|$.
Let $B(t) $ be an $m$-dimensional Brownian motion. If $A$ is a vector or matrix, its transpose is denoted by $A^T$. If $X \in \R^d$, then $|X|$ is the Euclidean norm. If $A$ is a matrix, its trace norm is denoted by $ |A | = \sqrt{ (A^T A)}$. For two real numbers $a$ and $b$,  $a \ve b:  = \max (a,b)$ and $a \we b : = \min (a,b)$. For a set $G$, its indicator function is denoted by $\II_G$. The scalar product of two vectors $X ,Y \in \R^d$ is denoted by $\lan X, Y \ran  $  or $X^T Y$. The largest integer which is less or equal to a real number $a$ is denoted by $ \lfloor  a \rfloor $.
In addition, we use $C$  to denote the generic constant that may change  from place to place.
\par
Let $\d: [0, +\infty ) \to [0, \tau]$ be the delay function which is Borel measurable.
Consider the following  SDDE of the form
\begin{align}\label{eq0}
dX(t) = f(X(t)),X(t-\d (t)))dt + g(X(t), X(t-\d (t)))dB(t), \; t \ge 0,
\end{align}
with the initial data
\begin{align}\label{eq3}
\{ X(\th) : -\tau \le \th \le 0  \}=\xi \in C([- \tau, 0]; \R^d),
\end{align}
where    $f: \R^d \K \R^d \to \R^d$ and $g : \R^d \K \R^d \to \R^{d\K m} $
are Borel-measurable functions.
We introduce  the following assumptions:
\begin{assumption}[Local Lipschitz Condition]\label{AA1_Local_Lips}
For any $R >0$, there exists a constant $L_R$ depending on $R$   such that
\begin{align*}
 | f(x_1,y_1) - f(x_2 ,y_2 )|^2 \ve  | g(x_1,y_1) - f(x_2 ,y_2 )|^2 \le L_R (|x_1 - x_2 |^2 + |y_1 - y_2 |^2),
\end{align*}
 for any $  x_1, x_2, y_1, y_2 \in \R^d$ with  $|x_1| \ve |x_2| \ve | y_2| \ve |y_2| \le R$.
\end{assumption}
\begin{assumption}[Generalized Khasminskii-type Condition]\label{AA3_Kha_cond}
There exist constants $K_1>0$, $K_2 \ge 0, K_3 \ge 0$ and  $ \b >2$  such that
\begin{align}\label{Guo_cond}
2 \langle x,   f(x,y)\rangle  +  |g (x,y)|^2  \le K_1 ( 1+ |x|^2 + |y|^2) - K_2 |x|^\b + K_3 |y|^\b, \qu \f x,y \in \R ^d.
\end{align}
\end{assumption}
\begin{assumption}\label{A2_bar_delta}
Let $\d: [0, +\infty ) \to [0, \tau]$  is continuously differentiable and there is  $\hat \d \in [0,1) $ such that
 $\displaystyle \left | \frac{d \d (t)}{d t} \right |  \le \hat \d $, for any $t \ge 0$.
\end{assumption}
\begin{lemma} \label{Lem1_2moments_of_X}
Suppose that Assumption \ref{AA1_Local_Lips}, \ref{AA3_Kha_cond}  and \ref{A2_bar_delta}  hold with
$\displaystyle K_2 > \frac{K_3}{(1- \hat \d)}  \ge 0$.
Then
for any given initial data \eqref{eq3}, there is a unique global solution $X(t)$ to \eqref{eq0} on $t \in [ -\tau, + \infty)$. Moreover, the solution $X(t)$ has the property that
\begin{align} \label{eq4}
\sup_{-\tau \le t \le T} \E |X(t)|^2 < \infty.
\end{align}
\end{lemma}
The proof of this result can be found in \citet{Song2013DCDS}.

\section{The truncated EM scheme for SDDEs with variable delay}
We first choose a strictly increasing continuous functions $\mu: \R_+ \to \R_+$ such that $\mu (R)\to \infty$ as $R\to \infty$ and
\begin{align}\label{eq21}
 \sup_{|x| \ve |x| \le r} \frac{|f(x,y)|}{(1+|x|+|y|)} \ve \frac{|g(x,y)|}{ (1+|x| + |y|)}  \le \mu (r) , \; \f r  \ge 1.
\end{align}
Denote by $\mu^{-1}$ the inverse function of $\mu$ and we see that $\mu^{-1}: [ \mu (1), \infty) \to \R_+$
is a strictly increasing continuous function. We then choose a constant $\hat h \ge  1 \ve \mu (1)$
and a strictly decreasing function $\varphi:(0,1] \to [\mu (1), + \infty)$ such that
\begin{align}\label{eq22}
\lim_{\D \to 0} \varphi(\D) = \infty \qu \textrm{and} \qu \D^{1/4} \varphi(\D) \le \hat h, \; \f \D \in (0,1].
\end{align}
For a given step size $\D \in (0,1]$, let us define a truncation mapping $\pd: \R^d \to \{ x \in \R^d: |x| \le \mu^{-1}(h(\D ) ) \}$ by
\begin{align}\label{eq23}
\pd(x) = \left ( |x| \we \mu^{-1}(\varphi(\D)) \right )\frac{x}{|x|}, \; \f x \in \R^d,
\end{align}
when $x=0$, we set $\displaystyle   x/ |x| = 0$. That is,  $\pd $ maps to itself if $|x| \le \mu^{-1} (\varphi(\D))$ and to $  \mu^{-1}(\varphi(\D)) x / |x| $ if $|x| > \mu^{-1} (\varphi(\D))$.
Define the following truncated functions
\begin{align}\label{eq24}
\fd (x,y) = f(\pd (x), \pd (y)) \qu \textrm{and} \qu  \gd (x,y) = g(\pd (x), \pd (y)), \; \f x,y \in \R^d.
\end{align}
From \eqref{eq21}, \eqref{eq23} and \eqref{eq24}, we have
\begin{align}\label{PP2}
|\fd (x,y)| \ve & |\gd (x,y)|  \le \varphi(\D)(1 + |\pd (x)| +|\pd (y)|) \le \varphi(\D)(1 + |x| +|y |),\; \f x,y \in \R^d,
\end{align}
which means that the truncated coefficients $f_\D$ and $g_\D$ grow at most  linearly  for a fixed step size $\D$, but $f$ and $g$ may not.
The following lemma shows  that these truncated coefficients can conserve the generalized Khasminskii-type condition for any $\D \in (0,1]$, the  proof
 of the lemma is similar to that of  \citet[Lemma 3.2]{Fei2020CPAA} and  so is omitted.
\begin{lemma}\label{Lem2_Conserving_Kcondition}
Let Assumptions   hold with $ K_2 \ge  K_3 \ge 0$. Then for any $\D \in (0,1]$,
\begin{align}\label{PP3}
2 \lan x, \fd (x,y) \ran +  |\gd (x,y)|^2 \le \hat K_1 (1 + |x|^2 +|y|^2) - K_2 |\pd (x)|^{\b}+ K_3 |\pd (y)|^{\b},\; \f x,y \in \R^d,
\end{align}
where $\hat K_1 =  2K_1 \left (  1 \ve [1/ \mu^{-1}(\varphi(1))] \right )$.
\end{lemma}
\par
Assume that   step size $\D \in (0,1]$ is a fraction of $\tau$. Take  $\D = \tau /M $ for some sufficiently large  integer $M$. 
Define $t_k = k \D$  and $ \d_k = \lfloor  \d (k \D) /\D \rfloor $, for any integer  $k\ge 0$.
%
%
 Then the boundedness of $\d$ gives
\begin{align}\label{eq3_1}
0 \le \d_k \le  \tau  /\D  =  M.
\end{align}
Define $ \displaystyle \k (t)  :=   \lfloor  t /\D \rfloor \D$, for any $t \ge - \tau$.
 The discrete-time truncated  EM scheme is defined as follows:
\begin{align}\label{Discrete_EM}
& y_{k+1} =   y_{k}  + \fd ( y_{k},y_{k-\d _k}) \D + \gd ( y_{k},y_{k-\d _k})  \D B_k, \;
  k \ge 0, \no \\
& y_ k = \x (t_k), \;   k = -M , -M + 1,  \cdots , 0,
\end{align}
 where  $\D B_k = B( \tkk )- B(\tk)$.
Define the continuous-time step approximations  $Z_1 (t)$ and $Z_2(t)$  by
 \begin{align}\label{Step_EM}
Z_1 (t) & = \sum_{k= -M}^{\infty} y_k \II_{[\tk, \tkk)}(t), \; \f  \; t \ge -\tau, \\
Z_2 (t) & = \sum_{k= 0}^{\infty} y_{k - \d_k} \II_{[\tk, \tkk)}(t), \; \f t \ge 0.
 \end{align}
where $\II$ is the indicator function. Define the   continuous-time  approximation $Y_{\D}(t)$ on
$t \in [-\tau, \infty)$ by
\begin{align}\label{Continuous_EM}
Y_{\D}(t) &  =  \x (0)  + \intt \fd (Z_1(s), Z_2(s))  ds  + \intt \gd (Z_1(s), Z_2(s))  dB(s), \;  t >0, \no \\
\yd (t) & = \xi (t), \;  - \tau \le t \le 0.
\end{align}
Thus $\yd (t)$  is an \Ito process  on $t \ge 0$ with \Ito differential
\begin{align}\label{eq_Yt}
d Y_{\D}(t) &  =  \fd (Z_1(t), Z_2(t))  dt  +  \gd (Z_1(t), Z_2(t))  dB(t).
\end{align}
It is useful to know that for any $t \in [\tk, \tkk)$ with $k \ge 0$,
 \begin{align}\label{Relations}
Z_1 (t) = \yd (\tk) = y_k \qu \textrm{and} \qu Z_2(t) = \yd ( \tk - \lf \d (\tk)/ \D \rf \D)= y_{k - \d_k},
 \end{align}
%
as well as
\begin{align}\label{YZ_relation}
\yd (t) - Z_1(t) & = \int ^t _ { \tk } \TDfds ds + \int ^t _ {\tk} \TDgds dB(s),
\end{align}
which means that the $Y_\D (t)$ and $Z_1 (t)$ coincide with the discrete solution at the grid points.
\begin{remark}
By \eqref{PP2},
we see from \eqref{Discrete_EM} that for a given step size $\D \in (0, 1]$, any $ p \ge 2$ and any $k \ge 1$,
\begin{align} \label{temp_eq1}
\E |\yk|^p \le C_{p, \| \xi \|, \hat h, \D}.
\end{align}
Moreover, this and  \eqref{YZ_relation}  guarantee that for  a given step size $\D \in (0, 1]$ and any $ p \ge 2$,
 \begin{align}\label{temp_eq3}
\E |\yd (t)|^p < \infty, \; \f t \ge 0.
 \end{align}
However, we cannot conclude that this bound is independent of $\D$. As a result of this observation, we need not apply stopping time arguments in the proof  Lemma \ref{Lem4_2moments_of_Y}, due to the fact that for any $\D \in (0, 1]$ and any $ p \ge 2$,
\begin{align}
 \E \int_0^T | \gd ( Z_1 (s), Z_2 (s))|^p ds < \infty.
\end{align}
In fact, by \eqref{Discrete_EM}, we have
\begin{align}
y_1 = y_0 + \fd (y_0, y_{- \d_0}) \D +  \gd (y_0, y_{- \d_0}) \D B_0.
\end{align}
Thus,  for  a given step size $\D \in (0, 1]$  and any $ p \ge 2$, by \eqref{PP2}, we have
\begin{align}
 \E | y_1|^p & \le 3^{p-1}\Big ( | \x (0)|^p + \E | \fd (y_0, y_{- \d_0})|^p \D^p +  \E | \gd (y_0, y_{- \d_0})|^p  \D^{0.5p} \Big )
 \no \\
 & \le 3^{p} (h (\D))^p \| \x \|^p \D ^{0.5p} \le 3^{p} \| \x \|^p  \hat  h^ p  \D^{p/4}.
\end{align}
Thus, by induction, we can show that \eqref{temp_eq1} holds. Moreover, rewrite \eqref{YZ_relation} as
\begin{align}
\yd (t) = y_k + \fd (y_k, y_{k - \d_k}) (t - \tk )+  \gd (y_k, y_{k -  \d_k}) (B(t) - B(\tk) ) , \; \f t \in [t_k, t_{k+1}), k \ge 0.
\end{align}
Then, for  a given step size $\D \in (0, 1]$  and any $ p \ge 2$
\begin{align*}
\E |\yd (t)|^p
& \le c_p \Big (  \E |y_k|^p + \E | \fd (y_k, y_{k - \d_k}) |^p \D^p +  \E |\gd (y_k, y_{k -  \d_k})|^p \D ^{0.5p} \Big ) \no \\
& \le c_p (\varphi(\D))^{p} ( \E |y_k|^p + \E |y_{k - \d _k}|^p)\D ^{0.5p} \no \\
& \le C_{p, \| \xi \|, \hat h, \D},   \; \f t \in [t_k, t_{k+1}), \;k \ge 0,
\end{align*}
which gives \eqref{temp_eq3}
\end{remark}

\section{Strong convergence}
\subsection{Strong convergence (without order) at time $T$}
\begin{lemma}\label{Lem3}
For any $\D \in (0,1]$ and $p>0$ ,
\begin{align}\label{eq4_1}
\E \left [ |\yd (t) - Z_1(t) |^{p}  | \F_{\k (t)} \right ]\le C   \D^{p/2 }( \varphi(\D))^{p }   (1 +|Z_1(t))|^{p} + |Z_2(t))|^{p}, \; \f t \ge 0,
\end{align}
where $C$ is a positive constant independent  of $\D$.
\end{lemma}
\pr
For any $p \ge 2$, by \eqref{PP2}, we see from \eqref{YZ_relation} that for any $ t \ge 0$
\begin{align*}
\E \Big [ |Y_\D (t) - Z_1 (t) |^{p}  | \F_{\k (t) }\Big ]
& = \E \left  [ \Big |  \int_{\k (t)}^t \fd (Z_1 (s), Z_2 (s))  ds
+  \int_{\k (t)}^t \gd (Z_1 (s), Z_2 (s))  dB(s)  \Big  |^{p}
 \Big | \F_{\k (t) } \right  ]   \no \\
 & \le C \E \left [ | \fd (Z_1 (t), Z_2 (t))(t -\k (t))  |^{p}    | \F_{\k (t) }\right ]
 +  C \E \left [ | \gd (Z_1 (t), Z_2 (t))( B(t) -B(\k (t)) )  |^{p}    | \F_{\k (t) }\right ] \no \\
 & \le C \D^{p} (\varphi(\D))^{p}  (1 + |Z_1 (t) |^{p} + |Z_2 (t) |^{p})
 + C  \D^{p /2 } (\varphi(\D))^{p  } (1 + |Z_1 (t) |^{p} + |Z_2 (t) |^{p}) \no \\
 & \le C  \D^{p /2 } (\varphi(\D))^{p  } (1 + |Z_1 (t) |^{p} + |Z_2 (t) |^{p}) ,
\end{align*}
this also holds for any $ p \in (0,  2) $  due to the \Holder inequality. $\Box$

\begin{lemma}\label{Lem_M}
Let Assumption \ref{A2_bar_delta} hold.
For any    $ k \in \{ 0,1,2, \cdots \} $, let $ k - \lf \d (k \D)/ \D \rf = u$, where  $u \in \{ -M, -M+1, \cdots, 0,1, \cdots, k\}$. Then
  \begin{align}\label{nt2}
\# \Big \{ j \in \{ 0,1,2, \cdots \}: j - \lfloor \d (j\D) /\D \rfloor  = u  \Big \}
\le \lfloor  (1- \hat \d)^{-1} \rfloor +1,
  \end{align}
where
$\# S$ denotes the number of elements of the set $S$.
\end{lemma}
This lemma provides an upper bound for the maximal number of indices $k \in \{1,2, \cdots  \}$
for which the expressions $k-\d _k$ are equal, the proof can be found in  \citet[Lemma 3]{MM2013}.
\begin{lemma}\label{Lem4_2moments_of_Y}
Let Assumptions \ref{AA1_Local_Lips},  \ref{AA3_Kha_cond} and \ref{A2_bar_delta} hold
with
$K_2 \ge ( \lfloor  (1- \hat \d)^{-1} \rfloor +1 ) K_3  \ge 0 $.
%
Then for any $\D \in (0,1]$,
\begin{align}\label{gs0}
 \sup_{ 0 \le t \le T} \E |Z_1 (t)|^2 \le C \qu \textrm{or} \qu  \sup_{ 0 \le k \D \le T} \E |\yk |^2 \le C,   ~~ \f T >0.
\end{align}
where $C$ is a positive constant independent  of $\D$.
\end{lemma}
\pr
For any integer $k \ge 0$, we conclude
from  \eqref{Discrete_EM} that
\begin{align}\label{eq11_1}
|y_{k+1}|^2 = |y_k|^2+ 2 \lan y_k, \fd(y_k, y_{k-\d_k}) \ran \D +
|\gd(y_k, y_{k-\d_k})|^2 \D  + | \fd(y_k, y_{k-\d_k})|^2 \D^2 + J_k,
\end{align}
where
\begin{align*}
J_k = 2 \lan y_k, \gd(y_k, y_{k-\d_k}) \D B_k \ran + 2 \lan \fd(y_k, y_{k-\d_k}), \gd(y_k, y_{k-\d_k}) \D B_k \ran\D
 + |\gd(y_k, y_{k-\d_k})|^2(|\D B_k|^2 -\D ).
\end{align*}
Obviously, $\E J_k =0$.
By \eqref{Lem2_Conserving_Kcondition}, we have
\begin{align}\label{gs1}
\E |y_{k+1}|^2 & \le \E |y_k|^2 + \hat K_1 \E (1 + |y_k|^2 + |y_{k- \d _k}|^2)\D + \E | \fd(y_k, y_{k-\d_k})|^2 \D^2 \no \\
& \qu  + \E  \Big [  - K_2 |\pd (y_k)|^\b  + K_3 | \pd (y_{k - \d _k})|^\b  \Big ]    \D , ~\f k \ge 0.
\end{align}
Moreover, by \eqref{PP2}, we have
\begin{align}\label{gs3}
& \E  | \fd( y_{k },y_{(k- \d _k)} ) |^2    \D^2
\le   (\varphi(\D))^2  \E (  1 + |y_k|^2 + |y_{k- \d _k}|^2 ) \D^2 \no \\
& \le \hat h ^2    \E ( 1 + |y_k|^2 + |y_{k- \d _k}|^2 ) \D^{3/2}
 \le \hat h ^2   \E ( 1 + |y_k|^2 + |y_{k- \d _k}|^2 )  \D .
\end{align}
Inserting \eqref{gs3} into \eqref{gs1} gives
\begin{align}\label{gs4}
\E |y_{k+1}|^2 & \le \E |y_k|^2 + ( \hat K_1 + \hat h^2 ) \E (1 + |y_k|^2 + |y_{k- \d _k}|^2)\D \no \\
& \qu  + \E \Big [  - K_2 |\pd (y_k)|^\b  + K_3 | \pd (y_{k - \d _k})|^\b    \Big ] \D , ~\f k \ge 0.
\end{align}
Thus, we  have
\begin{align}\label{gs5}
\E |y_{k}|^2 & \le  \| \x \|^2 + ( \hat K_1 + \hat h^2 )\sum_{j=0}^{k-1}\E (1 + |y_j|^2 + |y_{j- \d _j}|^2)\D \no \\
& \qu +
\E  \left [ \sum_{j=0}^{k-1} \Big (  - K_2 |\pd (y_j)|^\b  + K_3 | \pd (y_{j - \d _j})|^\b    \Big ) \right ] \D , ~\f k \ge 1.
\end{align}
By Lemma \ref{Lem_M}, we yields that
\begin{align}\label{eq5_2}
& \sum_{i=0}^{k-1}|\pd (y_{i-\d_i})|^{\b}  \D
 \le ( \lfloor  (1- \hat \d)^{-1} \rfloor +1 ) \sum_{j=-M}^{k-1}|\pd (y_{j})|^{\b} \D
 \no \\
 & = ( \lfloor  (1- \hat \d)^{-1} \rfloor +1 )  \sum_{j=-M}^{-1}|\pd ( y_{j}) |^{\b} \D + ( \lfloor  (1- \hat \d)^{-1} \rfloor +1 )  \sum_{j=0}^{k-1}|\pd (y_{j})|^{\b} \D  \no \\
 &
  \le ( \lfloor  (1- \hat \d)^{-1} \rfloor +1 ) \tau \| \xi \|^{\b} + ( \lfloor  (1- \hat \d)^{-1} \rfloor +1 )  \sum_{j=0}^{k-1}|\pd (y_{j})|^{\b} \D, \; \f k \ge 1,
\end{align}
where $\d_k = \lf \d (k\D) /\D \rf$.
Thus,
\begin{align}
& \sum_{j=0}^{k-1} \Big [ - K_2 |\pd (y_j)|^\b  + K_3 | \pd (y_{j - \d _j})|^\b    \Big ] \D \no \\
 & \le  K_3 ( \lfloor  (1- \hat \d)^{-1} \rfloor +1 )  \tau   \|\xi \|^{\b} -  ( K_2 -K_3 ( \lfloor  (1- \hat \d)^{-1} \rfloor +1 )    ) \sum_{j=0}^{k-1}|\pd (y_{j})|^{\b} \D  \no \\
 & \le K_3 ( \lfloor  (1- \hat \d)^{-1} \rfloor +1 )  \tau   \|\xi \|^{\b}.
\end{align}
Inserting this into \eqref{gs5}, we have
\begin{align}\label{gs6}
\E |y_{k}|^2
& \le ( \| \x \|^2 + K_3 ( \lfloor  (1- \hat \d)^{-1} \rfloor +1 )  \tau   \|\xi \|^{\b} ) + ( \hat K_1 + \hat h^2 )\sum_{j=0}^{k-1} \Big [ 1 + \E|y_j|^2 + \E|y_{j- \d _j}|^2  \Big  ] \D \no \\
& \le ( \| \x \|^2 + K_3 ( \lfloor  (1- \hat \d)^{-1} \rfloor +1 )  \tau   \|\xi \|^{\b} ) +
   ( \hat K_1 + \hat h^2 )\sum_{j=0}^{k-1} \Big [ 1 + 2 \sup_{-M \le i \le j}  \E|y_j|^2 \Big ] \D, ~ \f k \ge 1.
\end{align}
As this holds for any integer $k$ satisfying $1 \le k \le \lf T / \D \rf$,  while the sum of the right-hand-side (RHS) terms is non-decreasing
in $k$, we then have
\begin{align}\label{gs7}
\sup_{1 \le i \le k} \E |y_{i}|^2
& \le ( \| \x \|^2 + K_3 ( \lfloor  (1- \hat \d)^{-1} \rfloor +1 )  \tau   \|\xi \|^{\b} ) +
   ( \hat K_1 + \hat h^2 )\sum_{j=0}^{k-1} \Big [ 1 + 2 \sup_{-M \le i \le j}  \E|y_j|^2 \Big ] \D,
\end{align}
which implies that
\begin{align}
\sup_{-M  \le i \le k} \E |y_{i}|^2
& \le ( \| \x \|^2 + K_3 ( \lfloor  (1- \hat \d)^{-1} \rfloor +1 )  \tau   \|\xi \|^{\b} ) +
   ( \hat K_1 + \hat h^2 )\sum_{j=0}^{k-1} \Big [ 1 + 2 \sup_{-M \le i \le j}  \E|y_j|^2 \Big ] \D,\; \f k = 1,2, \cdots , \lf T / \D \rf.
\end{align}
By the discrete Gronwall inequality, we get  the desired assertion \eqref{gs0}. $\Box$
\begin{assumption} \label{AA4_initial_cond}
There is a pair of constants $K_4>0$ and  $\varrho \in (0,1]$ such that
\begin{align}
|\x (t) - \x (s)|   \le  K_4 |t-s|^{\varrho}, \; \forall s,t \in [-\tau ,0].
\end{align}
\end{assumption}
\begin{lemma}\label{Lem5}
Let Assumptions \ref{AA1_Local_Lips},  \ref{AA3_Kha_cond} and \ref{A2_bar_delta} hold
 with
 $K_2 \ge ( \lfloor  (1- \hat \d)^{-1} \rfloor +1 ) K_3  \ge 0 $.
 For any $R \ge \|  \xi \|$, define $\tau _R = \inf \{ t \ge 0: |X(t)|\ge R \}$
and
$\hat \rho _{\D, R} = \inf \{ t \ge 0: |Z_1 (t)|\ge R \}$.
Then
\begin{align}\label{asser21}
\PP ( \tau_R \le T) \le \frac{C}{R^2} \qu \textrm{and} \qu \PP ( \hat \rho_{\D,R} \le T) \le \frac{C}{R^2},\; \f T >0,
\end{align}
where $C$ is a positive constant independent  of $\D$.
\end{lemma}
\pr
By the \Ito formula and Assumption \ref{AA3_Kha_cond}, we have  that for any $t \in [0,T]$
\begin{align}\label{eq5_1_3}
\E |X (t \we \tau_R )|^2 & \le  |\x (0)|^2 + K_1 \E \int_0^{t \we \tau_R} \Big (1 + |X(s)|^2 + |X (s - \d (s))|^2 \Big ) ds  \no \\
& \qu  + \E \int_0^{t \we \tau_R} \Big ( -  K_2  |X(s)|^\b  +   K_3   |X(s - \d (s))|^\b  \Big ) ds  \no \\
& \le |\x (0)|^2  + K_1 T + K_1  \int_0^t \Big ( \E |X (s \we \tau_R)|^2 + \E | X( (s - \d (s)) \we \tau_R)|^2 \Big )ds   \no \\
& \qu +  \Big ( \frac{ K_3}{ 1- \hat \d} - K_2 \Big )\E \int_0^{t \we \tau_R} | X(s)|^\b ds + \frac{\tau \|  \x \|^\b }{ 1- \hat \d}  \no \\
& \le C + 2 K_1 \intt \Big (  \sup_{0 \le u \le s } \E | X (u \we \tau_R)|^2  \Big )ds,
\end{align}
where we have used the following estimates
\begin{align}\label{eq5_1_1}
& \int_0^{t \we \tau_R} |X (s - \d (s))|^\b  ds \le \frac{1}{1- \hat \d} \int_{ - \d (0)}^{ (t \we \tau_R ) - \d (t \we \tau_R )  } |X(u)|^\b du \no \\
&   \le \frac{1}{1- \hat \d} \int_{ - \tau}^{t \we \tau_R}|X(u)|^\b du
\le \frac{\tau \|  \x \|^\b }{ 1- \hat \d}  + \frac{1}{1- \hat \d} \int_{ 0}^{t \we \tau_R}|X(s)|^\b ds,
\end{align}
and
\begin{align}\label{eq5_1_11}
\intt \E |  X( ( s - \d (s) ) \we \tau_R)|^2 ds
& \le  \intt \Big (  \sup_{- \tau \le u \le s } \E | X (u \we \tau_R)|^2  \Big )ds  \no \\
& \le T  \| \xi \|^2 + \intt \Big (  \sup_{0 \le u \le s } \E | X (u \we \tau_R)|^2  \Big )ds.
\end{align}
Consequently,
$$\sup_{ 0 \le u \le t }\E |X (u  \we \tau_R )|^2
 \le C + 4 K_1 \intt \Big (  \sup_{0 \le u \le s } \E | X (u \we \tau_R)|^2  \Big )ds. $$
The Gronwall inequality gives
$$ \sup_{ 0 \le u \le t }\E |X (u  \we \tau_R )|^2  \le C.$$
Thus, $$ \E |X (T  \we \tau_R )|^2  \le C.$$
Finally, using the \Cheb inequality gives
$$ \displaystyle \PP ( \tau_R \le T) \le \frac{C}{R^2}. $$
Now, we begin to establish the second assertion in \eqref{asser21}.
The remaining  proof of this lemma is similar to that of \citet[Lemma 3.2]{LiMao2019IMA}, but more refined  techniques are  needed to overcome the difficilty
due to variable delay.
We observe that  $\hat \rho_{\D,R} = \vartheta_{\D, R}  \D $,
where $$\vartheta_{\D, R} : = \inf \{ k > 0 : |y_k| \ge R \}.$$
Clearly, $\hat \rho_{\D,R}$ and $\vartheta_{\D, R}$ are $\F_t $ and $\F_{\tk}$ stopping times, respectively.
It is useful to know  that
\begin{align}\label{Temp_eq2}
y_{(k+1) \we \vartheta_{\D, R}} - y_{k \we \vartheta_{\D, R}}  =  \II_{ \{ k <  \vartheta_{\D, R}  \}} (y_{k+1} - y_{k}) , \; \f k \ge 0,
\end{align}
see \citet[p. 477]{Prob_Book}.
Thus, from \eqref{Discrete_EM} and \eqref{Temp_eq2}, we  have
\begin{align*}
y_{(k+1)\we \vartheta_{\D, R} } = y_{k \we \vartheta_{\D, R}}  + \Big [\fd ( y_{k},y_{k-\d _k}) \D + \gd ( y_{k},y_{k-\d _k})  \D B_k  \Big ] \II_{ \{ k <  \vartheta_{\D, R}  \}}
\end{align*}
Consequently,
\begin{align}\label{eq5_3_3}
\E |y_{(k+1)\we \vartheta_{\D, R} }|^2 = & \E \Big [ |y_{k \we \vartheta_{\D, R} }|^2 +
 2 y_{k \we \vartheta_{\D, R} }^T \fd( y_{k \we \vartheta_{\D, R}},y_{(k- \d _k) \we \vartheta_{\D, R}} )\II_{ \{ k <  \vartheta_{\D, R}  \}}
 +  | \gd( y_{k \we \vartheta_{\D, R}},y_{(k- \d _k) \we \vartheta_{\D, R}} ) \D B_k  |^2   \II_{ \{ k <  \vartheta_{\D, R}  \}} \Big ] \no \\
 & + \E \Big [  | \fd( y_{k \we \vartheta_{\D, R}},y_{(k- \d _k) \we \vartheta_{\D, R}} ) |^2   \II_{ \{ k <  \vartheta_{\D, R}  \}}   \D^2  + \hat J_k  \Big ], \;  \f k \ge 0,
\end{align}
where
\begin{align*}
\hat J_k : & =  2 y_{k \we \vartheta_{\D, R} }^T  \gd( y_{k \we \vartheta_{\D, R}},y_{(k- \d _k) \we \vartheta_{\D, R}} ) \D B_k \II_{ \{ k <  \vartheta_{\D, R}  \}}   \no \\
& \qu + 2 \fd^T( y_{k \we \vartheta_{\D, R}},y_{(k- \d _k) \we \vartheta_{\D, R}} )  \gd( y_{k \we \vartheta_{\D, R}},y_{(k- \d _k) \we \vartheta_{\D, R}} ) \D B_k \II_{ \{ k <  \vartheta_{\D, R}  \}} \D.
\end{align*}
Note that
\begin{align}\label{eq5_3_4}
\D B_k \II_{ \{ k <  \vartheta_{\D, R}  \}} = B (t_{(k+1)\we \vartheta_{\D, R}} ) - B (t_{k \we \vartheta_{\D, R}} ).
\end{align}
Since $B(t)$ is a continuous  martingale, by the Doob martingale stopping time theorem, we have that
$ \E  \Big [  \D B_k \II_{ \{ k <  \vartheta_{\D, R}  \}} | \F_{t_{k \we \vartheta_{\D, R}}} \Big ] = 0$ and
for any $ A \in \R^{d \K m}$
\begin{align}\label{eq5_3_4_2}
\E  \Big [  | A \D B_k |^2 \II_{ \{ k <  \vartheta_{\D, R}  \}} | \F_{t_{k \we \vartheta_{\D, R}}} \Big ]
= |A|^2 \E  \Big [  ( t_{(k+1)\we \vartheta_{\D, R}}  - t_{k \we \vartheta_{\D, R}} ) | \F_{t_{k \we \vartheta_{\D, R}}} \Big ]
= |A|^2   \E \Big [ \II_{ \{ k <  \vartheta_{\D, R}  \}}   | \F_{t_{k \we \vartheta_{\D, R}}} \Big ] \D,
\end{align}
see \citet[p.12]{LiMao2019IMA}.
Then
\begin{align}\label{eq5_4_1}
\E \hat J_k &  =  2 \E \Big [ y_{k \we \vartheta_{\D, R} }^T  \gd( y_{k \we \vartheta_{\D, R}},y_{(k- \d _k) \we \vartheta_{\D, R}} )
 \E \big [  \D B_k \II_{ \{ k <  \vartheta_{\D, R}  \}} | \F_{t_{k \we \vartheta_{\D, R}}} \big ] \Big ] \no \\
 & \qu + 2 \E \Big [\fd^T( y_{k \we \vartheta_{\D, R}},y_{(k- \d _k) \we \vartheta_{\D, R}} )  \gd( y_{k \we \vartheta_{\D, R}},y_{(k- \d _k) \we \vartheta_{\D, R}} )
 \E \big [  \D B_k \II_{ \{ k <  \vartheta_{\D, R}  \}} | \F_{t_{k \we \vartheta_{\D, R}}} \big ] \Big ] \D
 \no \\
  & = 0.
\end{align}
Moreover, by \eqref{PP2} and Lemma \ref{Lem4_2moments_of_Y}, we have
\begin{align}\label{eq5_4_2}
& \E \Big [  | \fd( y_{k \we \vartheta_{\D, R}},y_{(k- \d _k) \we \vartheta_{\D, R}} ) |^2   \II_{ \{ k <  \vartheta_{\D, R}  \}} \Big ] \D^2
 = \E \Big [  | \fd( y_{k },y_{k- \d _k} ) |^2   \II_{ \{ k <  \vartheta_{\D, R}  \}} \Big ] \D^2 \no \\
& \le \E \Big [  | \fd( y_{k },y_{k- \d _k} ) |^2  \Big ] \D^2
\le  C (\varphi(\D))^2  \E \Big [  1 + |y_k|^2 + |y_{k- \d _k}|^2 \Big ]  \D^2 \no \\
& \le C \D^{3/2} \le C \D.
\end{align}
Plugging \eqref{eq5_4_1}, \eqref{eq5_4_2} into \eqref{eq5_3_3} and using Lemmas \ref{Lem2_Conserving_Kcondition}, \ref{Lem4_2moments_of_Y}, we have
\begin{align}\label{eq5_4_4}
\E |y_{(k+1)\we \vartheta_{\D, R} }|^2 &  \le  \E  |y_{k \we \vartheta_{\D, R} }|^2 + C \D
 + \E \Big [ \big ( 2 y_{k \we \vartheta_{\D, R}}^T \fd( y_{k \we \vartheta_{\D, R}},y_{(k- \d _k) \we \vartheta_{\D, R}} )
+  | \gd( y_{k \we \vartheta_{\D, R}},y_{(k- \d _k) \we \vartheta_{\D, R}} |^2  \big ) \II_{ \{ k <  \vartheta_{\D, R}  \}}  \Big ]\D  \no \\
&\le \E  |y_{k \we \vartheta_{\D, R} }|^2 + C \D +  \hat K_1 \E \Big [ \big ( 1 + |y_{k \we \vartheta_{\D, R}}|^2 + |y_{(k- \d _k) \we \vartheta_{\D, R}}|^2 \big ) \II_{ \{ k <  \vartheta_{\D, R}  \}}  \Big ] \D  + \hat H_k \no \\
& = \E  |y_{k \we \vartheta_{\D, R} }|^2 + C \D +  \hat K_1 \E \Big [ \big ( 1 + |y_{k }|^2 + |y_{(k- \d _k) }|^2 \big ) \II_{ \{ k <  \vartheta_{\D, R}  \}}  \Big ] \D  + \hat H_k \no \\
& \le \E  |y_{k \we \vartheta_{\D, R} }|^2 + C \D +   \hat K_1 \E \Big [ \big ( 1 + |y_{k }|^2 + |y_{k- \d _k }|^2 \big )  \Big ] \D  + \hat H_k \no \\
& \le \E  |y_{k \we \vartheta_{\D, R} }|^2 + C \D   + \hat H_k, \f k \ge 0,
\end{align}
where
\begin{align*}
\hat H_k = \E \Big [ \big ( - K_2 |\pd (y_k)|^\b  + K_3 | \pd (y_{k - \d _k})|^\b \big ) \II_{ \{ k <  \vartheta_{\D, R}  \}}   \Big ] \D.
\end{align*}
Thus, for any integer $k$ satisfying $  1 \le k \le \lf T/\D \rf  $,  we conclude from \eqref{eq5_4_4}  that

\begin{align}\label{eq5_5_1}
\E |y_{k \we \vartheta_{\D, R}}|^2 & \le  \| \x \|^2 + C  k \D + \sum_{j=0}^{k -1} \E \Big [ \big ( - K_2 |\pd (y_j)|^\b  + K_3 | \pd (y_{j - \d _j})|^\b \big ) \II_{ \{ j < \vartheta_{\D, R}  \}}   \Big ] \D    \no \\
& \le \| \x \|^2 + C T  +  \E \left [  \sum_{j=0}^{k -1}   \big ( - K_2 |\pd (y_j)|^\b  + K_3 | \pd (y_{j - \d _j})|^\b \big ) \II_{ \{ j < \vartheta_{\D, R}  \}}   \right ] \D  \no \\
&  =    \| \x \|^2 + C T  +  \E  \left  [ \sum_{j=0}^{ (k -1 ) \we ( \vartheta_{\D, R} -1)}   \big ( - K_2 |\pd (y_j)|^\b  + K_3 | \pd (y_{j - \d _j})|^\b \big ) \right  ] \D .
\end{align}
By Lemma \ref{Lem_M}, we get
\begin{align}
& \sum_{j = 0}^{ (k-1) \we ( \vartheta_{\D, R} - 1)} |\pd (y_{j - \d _j})|^\b
\le ( \lfloor  (1- \hat \d)^{-1} \rfloor +1 ) \sum_{i = - M}^{ (k-1) \we ( \vartheta_{\D, R} -1)} |\pd (y_{i })|^\b \no \\
& = ( \lfloor  (1- \hat \d)^{-1} \rfloor +1 ) \sum_{i = -M }^{-1}  |\pd (y_{i })|^\b +  ( \lfloor  (1- \hat \d)^{-1} \rfloor +1 ) \sum_{i = 0  }^{(k-1) \we ( \vartheta_{\D, R}-1) }  |\pd (y_{i })|^\b \no \\
& \le  ( \lfloor  (1- \hat \d)^{-1} \rfloor +1 ) M \| \x \|^\b  +  ( \lfloor  (1- \hat \d)^{-1} \rfloor +1 ) \sum_{i = 0  }^{(k-1) \we ( \vartheta_{\D, R}-1) }  |\pd (y_{i })|^\b, \; \f k \ge 1 .
\end{align}
Consequently,
\begin{align}\label{eq5_5_2}
& \sum_{j=0}^{ (k -1 ) \we ( \vartheta_{\D, R} -1 )}  \Big [ \big ( - K_2 |\pd (y_j)|^\b  + K_3 | \pd (y_{j - \d _j})|^\b \big ) \Big ] \D   \no \\
&   \le K_3 ( \lfloor  (1- \hat \d)^{-1} \rfloor +1 ) ( M  \D )\| \x \|^\b  +  ( K_3 ( \lfloor  (1- \hat \d)^{-1} \rfloor +1 ) - K_2 )  \sum_{j = 0}^{ (k-1) \we ( \vartheta_{\D, R} -1)} |\pd (y_{j})|^\b \D
  \no \\
  & \le K_3 ( \lfloor  (1- \hat \d)^{-1} \rfloor +1 ) \tau \| \x \|^\b,\; \f k \ge 1 .
\end{align}
Inserting this into \eqref{eq5_5_1} gives
\begin{align*}
\E |y_{k \we \vartheta_{\D, R}}|^2 & \le
  \| \x \|^2 + c_1 T  +  K_3 ( \lfloor  (1- \hat \d)^{-1} \rfloor +1 ) \tau \| \x \|^\b ,\; \f k = 1,2, \cdots, \lf T/\D \rf  .
\end{align*}
In particular, we have
$$ \E |y_{\lf T/\D \rf   \we \vartheta_{\D, R}}|^2  \le C, $$
or equivalently,
$$ \E |Z_1 (T \we \hat \rho_{\D,R} )|^2  \le C $$
which  implies that
$$ {R^2} \PP ( \hat \rho_{\D,R} \le T) \le \E \Big [ \II_{ \{ \hat \rho_{\D,R} \le T   \}}  |Z_1 (\hat \rho_{\D,R} )|^2   \Big ]  \le  \E |Z_1 (T \we \hat \rho_{\D,R} )|^2  \le C .  $$
Thus, the proof is finished. $\Box$
\par

\begin{remark}
It should be pointed out that if we use the usual continuous proof of \citet[Lemma 3.3]{Guo2018NA} to estimate the second assertion in \eqref{asser21},
then there will be a term $J^\star$ with the following form we have to address.
\begin{align}\label{Temp_eq1}
J^\star : = \int_0^{ \th} | \pd  ( Z_2 (s) )|^\b ds -  \bar \k \int_0^{ \th} | \pd  ( Z_1 (s) )|^\b ds,
\end{align}
where $\bar \k =  \lfloor  (1- \hat \d)^{-1} \rfloor +1$ and $\theta$ is a stopping time.
Of course, by the known conditions, we can show that
\begin{align}\label{Temp_eq2}
J^\star
& \le \bar \k  \tau \| \x \|^\b
+  \Big (  | \pd ( y_{\lf \th /\D \rf  - \d_{\lf \th /\D \rf}})|^\b   - \bar \k | \pd (y_{ \lf \th /\D \rf   })|^\b  \Big )  (\th - \lf \th /\D \rf).
\end{align}
But we see from \eqref{Temp_eq2} than this estimate remains a
"tail", namely, the second term on the  right hand side of \eqref{Temp_eq2}, and we have on other method to  address   $J^\star$  properly.
%
However,
if we take $\th$ to be the grid point, then
the "tail" vanishes. This motivates the above discrete proof in Lemma \ref{Lem5}. \\

%
\end{remark}
The following theorem establish the strong  convergence (without  order) results of the truncated EM method.
\begin{theorem}\label{Th1_convergnce}
Let Assumptions \ref{AA1_Local_Lips}, \ref{AA3_Kha_cond}, \ref{A2_bar_delta} and    \ref{AA4_initial_cond} hold
 with
 $K_2 \ge ( \lfloor  (1- \hat \d)^{-1} \rfloor +1 ) K_3  \ge 0 $.
Then for any $ q \in [1,2)$,
  \begin{align}\label{asser2}
\lim_{\D \to 0} \E | X (T) - Z_1(T)|^q = 0, \; \f T>0.
 \end{align}
\end{theorem}
\pr Let $\tau_R $ and $\hat \rho_{\D, R}$ be the same as before.  Let $q \in [1,2)$.
 Define $ \hat  \th_{\D, R} = \tau_R \we \hat \rho_{\D, R}$ and
$ \hat e_\D (T) = X(T) - Z_1 (T)$, for any $T >0$.
 By the Young inequality, for any $\eta >0$, we have
\begin{align*}
\E |\hat e_\D (T)|^q &  = \E  \big [ | \hat e_\D (T)|^q  \II_{ \{  \hat \th_{\D, R} >T+1 \} }  \big  ]
                    + \E  \big [ | \hat e_\D (T)|^q  \II_{ \{ \hat \th_{\D, R} \le T +1 \} }  \big  ]  \no \\
                    & \le \E  \big [ | \hat e_\D (T)|^q  \II_{ \{  \hat \th_{\D, R} >T+1 \} }  \big  ] + \frac{ q \eta}{2} \E |\hat e_\D (T)|^2
                   +  \frac{2-q}{2 \eta ^{ q/(2 - q)}} \PP ( \hat \th_{\D,R} \le T+1).
\end{align*}
In this theorem, $C_R$ denotes a positive constant depending  on $R$ but independent of $\D$, its value may be different for different appearance.
By Lemmas \ref{Lem1_2moments_of_X} and \ref{Lem4_2moments_of_Y}, we get that
$$ \E |\hat e_\D (T)|^2 \le 2 \E |X(T)|^2 + 2 \E |Z_1(T)|^2 \le C.  $$
While from Lemma \ref{Lem5}, we have
$$\PP ( \th_{\D,R} \le T+1) \le \PP ( \tau_R \le T+1) + \PP ( \hat \rho_{\D,R} \le T+1)   \le \frac{C}{R^2}. $$
Consequently, we have
\begin{align}\label{eq6_1_1}
\E |\hat e_\D (T)|^q & \le  \frac{ q \eta C }{2} + \frac{(2-q) C }{2 R^2\eta ^{ q/(2 - q)}} + \E  \big [ | \hat e_\D (T)|^q \II_{ \{  \hat \th_{\D, R} >T+1 \} }  \big  ].
\end{align}
Let $\hat \eps >0  $ be arbitrary. Choose $\eta > 0 $ sufficiently small for   $ \displaystyle   \frac{ q \eta C }{2} \le \hat \eps  $
and then choose $R$ sufficiently large for $ \displaystyle  \frac{(2-q) C }{2 R^2\eta ^{ q/(2 - q)}} \le \hat \eps$.
Then for such chosen $R$, we see from \eqref{eq6_1_1} that
\begin{align*}
\E |\hat e_\D (T)|^q \le \E  \big [ | \hat e_\D (T)|^q  \II_{ \{  \hat \th_{\D, R} >T+1 \} }  \big  ] + 2 \hat \eps.
\end{align*}
If we can show that
\begin{align}\label{eq5_7_1}
 \lim_{ \D \to 0}  \E  \big [ | \hat e_\D (T)|^q  \II_{ \{  \hat \th_{\D, R} >T+1 \} }  \big  ] =0,
\end{align}
the desired assertion \eqref{asser2} follows.  Define the truncated functions
$$ F_R (x,y) = f \Big   ( (|x| \we R) \frac{x}{|x|},  (|y| \we R) \frac{y}{|y|}\Big  )  $$
 and
$$ G_R (x,y) = g \Big   ( (|x| \we R) \frac{x}{|x|},  (|y| \we R) \frac{y}{|y|}\Big  ), $$
for any $x,y \in \R^d$.
Without loss of any generality, we assume that $\D^* $ is sufficiently small for
$ \mu^{-1} ( h(\D^*))  \ge R$.  Then, for any $\D \in (0, \D^*]$, we get that
\begin{align*}
\fd (x,y) = F_R (x,y) \qu \textrm{and} \qu \gd (x,y) = G_R (x,y)
\end{align*}
for any $x, y \in \R^d$ with $|x| \ve |y| \le R$. Consider the following SDDE
\begin{align}\label{eq5_7_2}
d z(t) = F_R (z(t) , z(t - \d (t))) dt +G_R (z(t) , z(t - \d (t)))dB(t), \; t \ge 0,
\end{align}
with the initial data $z(t) = \x (t)$ on $t \in [- \tau  ,0]$.
By Assumption \ref{AA1_Local_Lips}, we observe that $F_R(x,y)$ and $G_R(x,y)$  are globally Lipschitz  continuous with
the Lipschitz constant  $ L_R$.  Hence, SDDE \eqref{eq5_7_1} has a unique global solution $z(t)$ on $t \ge - \tau$ satisfying
\begin{align}\label{eq5_7_3}
\PP (  z(t \we \tau_R) = X (t \we \tau_R)  \; \textrm{for} \; \textrm{any} \;t \in [0,  T]   ) = 1.
\end{align}
On the other hand, for any $\D \in (0, \D^* ]$, we apply the (classical) EM method to the SDDE \eqref{eq5_7_2} and
 we denote  $z_\D (t)$ and $\bar z_\D (t)$  by the continuous-time continuous-sample and  the piecewise constant EM solutions, respectively.  Then we see from  \citet[Theorem 2.1]{Mao2003SDDE}   that continuous-time continuous-sample EM solution $z_\D (t)$ has the property
\begin{align}\label{eq5_7_4}
\E \Big [ \sup_{ 0 \le t \le T} | z(t) - z_\D (t)|^q   \Big ] \le c_2 \D^{q (0.5 \we \varrho)},
\end{align}
where $c_2$ is a positive constant dependent of $ L_R$, $T$, $\x$, $q$ but independent of $\D$.
From this and the fact that
\begin{align}\label{eeq1}
\E |z_\D (T \we ( \hat \rho_{\D, R} -1)   ) - \bar z_\D (T \we ( \hat \rho_{\D, R} -1) )| ^q \le C_R \D ^{0.5(q \we \varrho ) },
\end{align}
see \citet[Corollary 3.4]{Mao2003SDDE}, we conclude that
\begin{align}\label{eeq2}
\E |z(T \we ( \hat \rho_{\D, R} -1)) - \bar z_\D (T \we ( \hat \rho_{\D, R} -1))| ^q \le C_R \D ^{0.5(q \we \varrho )}.
\end{align}
Moreover,
\begin{align}\label{eq5_7_5}
\PP \Big (  Z_1(t \we ( \hat \rho_{\D, R} -1) ) = \bar z_\D (t \we ( \hat \rho_{\D, R} -1 ) )  \; \textrm{for} \; \textrm{any} \; t \in [0,T]  \Big ) = 1.
\end{align}
Consequently,
 \begin{align}
 & \E  \big [ | \hat e_\D (T)|^q  \II_{ \{  \hat \th_{\D, R} >T+1 \} }  \big  ]
 =  \E  \big [ | \hat e_\D (T \we ( \hat \th_{\D, R} -1 )) |^q  \II_{ \{  \hat \th_{\D, R}-1 >T \}  } \big  ] \no \\
 & \le \E  \big [ | X  (T \we ( \hat \th_{\D, R}  -1 )) - Z_1  (T \we ( \hat \th_{\D, R} -1 )) |^q   \big  ] \no \\
 & = \E  \big [ | z  (T \we ( \hat \th_{\D, R} -1)) - \bar z_\D  (T \we ( \hat \th_{\D, R} -1 )) |^q   \big  ]  \no \\
  & \le  C_R \D^{q (0.5 \we \varrho)},
 \end{align}
which establish   \eqref{asser2}.
  Thus, the proof is finished.
   $\Box$
\subsection{Strong convergence order at time $T$}

In this section, we mainly discuss the strong  convergence order of the truncated EM method for \eqref{eq0}. We need a stronger condition than \eqref{Guo_cond}.
\begin{assumption} \label{AA3_Kha_cond_P0}
There is a pair of constants $p_0 >2$ and $\bar K_1  >0$ such that
\begin{align*}
 \lan x, f(x,y) \ran + \frac{p_0 -1}{2} |g(x,y)|^2 \le \bar K_1 (1 + |x|^2 + |y|^2) , \; \f x_1,x_2 \in \R^d.
\end{align*}
\end{assumption}
We cite a known result as a Lemma, see e.g., \cite{Mao2005SAA}.
\begin{lemma} \label{Lem1_2moments_of_P0}
Suppose that Assumption \ref{AA1_Local_Lips}  and \ref{AA3_Kha_cond_P0} hold. Then
for any given initial data \eqref{eq3}, there is a unique global solution $X(t)$ to \eqref{eq0} on $t \in [-\tau, + \infty)$. Moreover, the solution $X(t)$ has the property that
\begin{align*}
\sup_{-\tau \le t \le T} \E |X(t)|^{p_0} < \infty.
\end{align*}
\end{lemma}
\begin{lemma}\label{Lem4.2_Conserving_Kcondition}
Let Assumption \ref{AA3_Kha_cond_P0} hold. Then for any $\D \in (0,1]$,
\begin{align*}
\lan x, \fd (x,y) \ran + \frac{p_0 -1 }{2} |\gd (x,y)|^2 \le 2 \bar K_1 \left (  1 \ve [1/ \mu^{-1}(\varphi(1))] \right )(1 + |x|^2 +|y|^2) , \; \f x_1,x_2 \in \R^d.
\end{align*}
\end{lemma}
The proof can be found in \citet[Lemma 3.2]{Fei2020CPAA}.
\par
Using the similar techniques in the proofs of Lemma \eqref{Lem4_2moments_of_Y}, we obtain the following lemma which provides an upper bound for the $p_0$-th moment of the numerical solution $Y_\D$.
\begin{lemma}\label{Lem4_4p_0_th moments_of_Y}
Let Assumptions \ref{AA1_Local_Lips}, \ref{AA3_Kha_cond_P0}, \ref{A2_bar_delta} hold. Then
\begin{align*}
\sup_{0 <  \D \le 1} \sup_{ 0 \le t \le T} \E |\yd (t)|^{p_0} \le C,
\end{align*}
where $C$ is a positive constant independent of $\D$.
\end{lemma}
%
Let $\mathcal U $ be the family of continuous function $ U:\R^d \K \R^d \to \R_+$ such that
for any $b>0$, there exists a positive constant $\k_b$ for which
$$ U(x_1,x_2) \le \k_b |x_1 - x_2|^2, \; \f x_1,x_2 \in \R^d.$$
\begin{assumption}[Global Monotonicity with $U$ Function and Polynomial Growth conditions]\label{AS31_monotonity_and_Poly_cond}
There exist constants $p_1 >2$, $l \ge 0$ and $ \bar K_2 >0$ as well as a function $u \in \mathcal U$ such that
\begin{align}\label{eq16_1}
& \lan x_1 - x_2, f(x_1,y_1) - f(x_2,y_2) \ran + \frac{p_1 -1}{ 2} |g(x_1,y_1) - g(x_2,y_2)|^2  \no \\
& \le \bar K_2 (|x_1-x_2|^2 + |y_1-y_2|^2) - \frac{1}{1- \hat \d} U(x_1,x_2) + U(y_1,y_2),\; \f x_1,y_1,x_2,y_2 \in \R^d
\end{align}
and
\begin{align}\label{eq16_2}
&|f(x_1,y_1) - f(x_2,y_2)|^2 \ve |g(x_1,y_1) - g(x_2,y_2)|^2  \no \\
& \le \bar K_2(1+|x_1|^l + |x_2|^l + |y_1|^l + |y_2|^l)   (|x_1-x_2|^2 + |y_1-y_2|^2), \; \f x_1,y_1,x_2,y_2 \in \R^d.
\end{align}
\end{assumption}
Note that
\begin{align}\label{temp1017}
|\pi_\D (x)| \le |x|, \qu |\pi_\D (y)|\le |y|, \qu |\pd (x) - \pd (y)|^2 \le |x-y|^2, \; \f x,y \in \R^d.
\end{align}
Thus, from  \eqref{eq16_2}, we have the following growth condition
\begin{align}\label{eq16_4}
|f(x,y)|^2 \ve |g(x,y)|^2 \le \bar K_3 (1 + |x|^{2+l} + |y|^{2+l}), \;\f x,y \in \R^d,
\end{align}
and
\begin{align}\label{eq16_5}
& |\fd(x,y)|^2  = |f(\pd (x), \pd(y))|^2 \no \\
& \le \bar K_3 (1 + |\pd (x)|^{2+l} + |\pd (y)|^{2+l})
 \le \bar K_3 (1 + |x|^{2+l} + |y|^{2+l}), \;\f x,y \in \R^d,
\end{align}
as well as
\begin{align}\label{eq16_6}
|\gd(x,y)|^2
& \le \bar K_3 (1 + |x|^{2+l} + |y|^{2+l}), \;\f x,y \in \R^d,
\end{align}
where $\bar K_3$ is a positive constant depending  on $\bar K_2$.
\begin{lemma}\label{Lem3.4}
Let Assumptions \ref{AA1_Local_Lips}, \ref{AA3_Kha_cond_P0}, \ref{A2_bar_delta}, \ref{AA4_initial_cond} hold with $ p_0 \ge 2 +  l$. Then  for any $\D \in (0,1]$ and  any $\displaystyle  p \in \left [2, \frac{p_0}{ 1 + l/2} \right ]$,
\begin{align}\label{eq19_2}
\E | \yd (t) - Z_1(t)|^p \le C \D^{0.5p}, \; \f t \ge 0,
 \end{align}
and
\begin{align}\label{eq19_3}
\E | \yd (t- \d (t) ) - Z_2(t)|^p \le C \D^{( \varrho \we 0.5 )p}, \; \f t \ge 0,
 \end{align}
where $C$ is a positive constant independent of $\D$.
\end{lemma}
\pr
Let $\displaystyle p \in \left [ 2, \frac{p_0}{ 1 + l/2} \right ]$.
For any  $ t \in [t_k, t_{k+1})$ with $k \ge 0$,
we see from \eqref{YZ_relation} that
\begin{align*}
\E |\yd (t) - Z_1 (t)|^p & = \E |\yd (t) - \yd (\tk)|^p
\no \\
& \le C \left (  \D^p \E |\fd (\yk, y_{k- \d _k})|^p+
\D^{0.5 p} \E |\gd (\yk, y_{k- \d _k})|^p
\right ).
\end{align*}
By \eqref{eq16_5}, we have
\begin{align}\label{eq19_5}
\E |\fd (\yk, y_{k- \d _k})|^p & \le C \E \left ( 1 + |\yk|^{2+l} + |y_{k-\d_k}|^{2 + l} \right )^{p/2} \no \\
& \le C \left ( 1 + \E |\yk|^{(1+ l/2)p} + \E |y_{k-\d_k}|^{(1+ l/2)p} \right ) \le C,
\end{align}
where  Lemma \ref{Lem1_2moments_of_P0} has been used.
Similarly, we can show that
\begin{align}\label{eq19_6}
\E |\gd (\yk, y_{k- \d _k})|^p \le C.
\end{align}
Thus, $$\E |\yd (t) - Z_1 (t)|^p \le C\D^p + C\D^{0.5p} \le C \D^{0.5p}.$$
We now begin to establish assertion \eqref{eq19_3}.
Recall that  $\d_k = \lf \d (\tk)/\D \rf$ and
\begin{align}\label{eq21_1}
\yd (t - \d (t)) - Z_2 (t) = \yd (t - \d (t)) - \yd ( (k - \d_k)\D).
\end{align}
By Assumption \ref{A2_bar_delta}, we have the following useful estimate
\begin{align}\label{eq21_3}
|(t - \d (t)) - (k - \d_k)\D | \le ( \lf \hat \d \rf + 4) \D,
\end{align}
see \citet[Lemma 4.6]{Deng2019SSCE}.
From now on, we will use  $c_p$  to denote  genetic positive constants dependent only on $p$   and its values may change between occurrences.
Now consider the following four possible cases.
\begin{enumerate}
   \item[\textbf{Case 1:}] If $t - \d (t) \ge (k - \d_k)\D \ge 0$  or $(k - \d_k)\D  \ge t - \d (t) \ge 0$, then
   it follows from \eqref{eq21_1} that
   \begin{align}\label{eq21_4}
\E |\yd (t- \d(t)) - Z_2 (t)|^p & = \E \left | \int_{(k-\d_k)\D}^{t - \d (t) } \TDfds ds  + \int_{(k-\d_k)\D}^{t - \d (t) } \TDgds dB(s) \right |^p \no \\
& \le  c_p |t-\d(t) - (k-\d_k)\D |^{p-1} \int_{(k-\d_k)\D}^{t - \d (t) } \E |\TDfds |^p ds \no \\
& \qu +  c_p |t-\d(t) - (k-\d_k)\D |^{0.5p -1} \int_{(k-\d_k)\D}^{t - \d (t) } \E |\TDgds |^p ds \no \\
& \le C |t-\d(t) - (k-\d_k)\D |^{0.5p }
& \no \\
& \le  C  (\lf \hat \d \rf+4)^{0.5p} \D^{0.5p},
   \end{align}
   where  the \BDG inequality, \eqref{eq19_5}  \eqref{eq19_6} and  \eqref{eq21_3} have been used.
   \item[\textbf{Case 2:}] If $t - \d (t)  \le (k-\d_k)\D  \le 0 $ or $ (k-\d_k)\D \le  t - \d (t)    \le 0 $, by Assumption \ref{AA4_initial_cond} and \eqref{eq21_3} we have
   \begin{align}\label{eq22_1}
\E |\yd (t- \d(t)) - Z_2 (t)|^p & = |\x (t - \d (t) ) - \x ( (k-\d_k)\D) |^p \no \\
& \le \left ( K_4^p (\lf \hat \d \rf+4)^{\varrho p } \right )\D^{\varrho p} .
   \end{align}
   \item[\textbf{Case 3:}] If $t - \d (t) \ge 0 \ge (k-\d_k)\D  $, then
   \begin{align*}
 t - \d (t) \le  (\lf \hat \d \rf+4)\D  \qu \textrm{and}\qu
  -(k-\d_k)\D \le (\lf \hat \d \rf+4)\D .
   \end{align*}
     Thus, we have
  \begin{align}\label{eq22_4}
\E |\yd (t- \d(t)) - Z_2 (t)|^p & = \E |\yd (t- \d(t)) - \yd ((k-\d_k)\D)|^p  \no \\
& \le 2^{p-1} \E | \yd (t - \d(t))- \xi (0) |^p + 2^{p-1} \E | \x (0) - \x ( (k - \d _k)\D ) |^p \no \\
& \le C (t - \d(t))^{0.5p} + C K_4^p (- (k-\d(k))\D)^{\varrho p} \no \\
& \le C  \left ( (\lf \hat \d \rf+4)^{\varrho p } + (\lf \hat \d \rf+4)^{0.5 p } \right )  \D^{ (\varrho \we 0.5) p}.
  \end{align}
   \item[\textbf{Case 4:}] If $  (k-\d_k)\D \ge 0 \ge   t - \d (t) $, in a similar way as  \eqref{eq22_4} was obtained,
   we also have
     \begin{align}\label{eq22_5}
\E |\yd (t- \d(t)) - Z_2 (t)|^p
& \le C \left (   (\lf \hat \d \rf+4)^{0.5 p } + (\lf \hat \d \rf+4)^{\varrho p } \right )  \D^{ (\varrho \we 0.5) p}.
\end{align}
\end{enumerate}
Combining these different cases together,
we obtain the desired assertion \eqref{eq19_3}. $\Box$

We then can obtain the optimal rate of strong convergence for truncated EM approximation.
\begin{theorem}\label{Th_3__convergence_rate} 
Let Assumptions  \ref{AA3_Kha_cond_P0}, \ref{AS31_monotonity_and_Poly_cond}, \ref{A2_bar_delta}  and \ref{AA4_initial_cond} hold  with $ p_0 \ge 2 + 3 l$. Then for any $\D \in (0,1]$,
\begin{align}\label{eq17_5}
 \E |X(T) - \yd (T)|^2 \le C \left ( \D^{2 \varrho \we 1 }  \ve [\mu^{-1} (\varphi (\Delta))]^{-(p_0 - l -2)}\right ), \; \f T>0,
\end{align}
and
\begin{align}\label{eq17_6}
 \E |X(T) - Z_1 (T)|^2 \le C \left ( \D^{2 \varrho  \we 1}  \ve [\mu^{-1} (\varphi (\Delta))]^{-(p_0 - l -2)}\right ), \; \f T>0,
\end{align}
where $C$ is a positive constant independent of $\D$.
In particular, let
\begin{align}\label{phy11}
\mu (r) = \bar K_3 ^{1/2}r^{l/2}, \; \f r \ge 1 \qu \textrm{and}\qu \varphi(\D) = \hat h \D^{-1/4}, \; \f \D \in(0,1].
\end{align}
%
Then for any $\D\in (0,1]$,
\begin{align}\label{asser31}
 \E |X(T) - \yd (T)|^2  \le C \D^{1 \we 2 \varrho  }\qu \textrm{and} \qu  \E |X(T) - Z_1(T)|^2  \le C \D^{1 \we 2 \varrho }, \; \f T>0.
\end{align}

\end{theorem}
\begin{remark}\label{rem1}
If  constant delay is considered in \eqref{eq0},  then Theorem \ref{Th_3__convergence_rate} reduces to \citet[Theorem 3.6]{Fei2020CPAA}.
Comparing with the two theorems, we observe  the following differences:
\begin{itemize}
  \item Theorem \ref{Th_3__convergence_rate} requires a slightly stronger condition on the parameters, namely, $p_0 \ge 2 +3l$, while
  \citet[Theorem 3.6]{Fei2020CPAA} requires $p_0 > 2 +l$;
  \item Theorem \ref{Th_3__convergence_rate} allows a slightly weaker control function $\mu$ for truncation, namely, $\mu(r)=Cr^{l/2}$, while
  \citet[Theorem 3.6]{Fei2020CPAA} allows $\mu(r) = Cr^{1 + l/2}$;
  \item In Theorem \ref{Th_3__convergence_rate} the truncated EM approximation $\yd (T)$ achieves a  better order of  $\mathcal L^2$-convergence which is  $2 \varrho  \we 1  $, while in  \citet[Theorem 3.6]{Fei2020CPAA}  the corresponding convergence order is  $2 \varrho  \we (1-2\varepsilon) \we 2 \varepsilon (p_0 - l -2)/(2+l)  $ for some $\varepsilon \in (0,1/4]$, which can be  close to the rate of Theorem \ref{Th_3__convergence_rate}, but  $p_0$ should be required to  sufficiently large.
%
\end{itemize}
\end{remark}
\begin{remark}\label{rem2}
It is worth mentioning how our work compares with that of \citet{GuoZhan2018IJCM}, who proved the strong convergence results of
partially truncated EM method applied to  the SDDEs with variable delay and Markovian switching. What differentiates  our work from \cite{GuoZhan2018IJCM} are:
\begin{itemize}
  \item We remove the restrictive condition (3.25) in \cite[Theorem 3.12]{GuoZhan2018IJCM}, i.e.,
  \begin{align}\label{guo_cond4}
h(\D) \ge \mu \Big ( (\D^{2\nu} \ve \D \varphi^2 (\D) )^{-2/(p-2)}   \Big ),
  \end{align}
  which could force the step size to be so small that the truncated EM method would be inapplicable;
   \item We relax the grown constraint on the Khasminskii-type condition (2.13) in \cite[Assumption 2.3]{GuoZhan2018IJCM}, i.e.,
   \begin{align}\label{guo_cond5}
2 \langle x,   f(x,y)\rangle  +  |g (x,y)|^2  \le K ( 1+ |x|^2 + |y|^2), \qu \f x,y \in \R ^d,
   \end{align}
 by a generalized  Khasminskii-type condition \eqref{Guo_cond}, and
   global monotonicity condition (3.15) in  \cite[Assumption 3.9]{GuoZhan2018IJCM}, i.e.,
\begin{align}\label{guo_cond6}
& \lan x_1 - x_2, f(x_1,y_1) - f(x_2,y_2) \ran + \frac{p_1 -1}{ 2} |g(x_1,y_1) - g(x_2,y_2)|^2  \no \\
& \le  K (|x_1-x_2|^2 + |y_1-y_2|^2) ,\; \f x_1,y_1,x_2,y_2 \in \R^d
\end{align}
   by the global monotonicity condition with $U$ function \eqref{eq16_1};
   \item We obtain the optimal  order $1/2$ of strong convergence which is higher than that of \citet[Theorem 3.12]{GuoZhan2018IJCM}, i.e.,
   \begin{align*}
\E |X(T) - Y_\D (T)|^2 \le C \Big ( \D \ve \D \varphi^2 (\D) \Big),
   \end{align*}
  under a slightly stronger condition on the parameters.

\end{itemize}

\end{remark}
\textbf{Proof of Theorem  \ref{Th_3__convergence_rate}.}
Let  $ p_0 \ge 2 + 3 l$ and $R \ge \| \x \|$, define
the stopping time $$ \rho _R = \inf \{ t \ge 0 : |X(t)| \we |\yd (t)| \ge R  \}   .$$  Set $e_\D (t) = X(t) - \yd (t)$, for any $t \in [-\tau, T]$, which means  that $e_\D (t) = 0$, for any $t \in [-\tau, 0]$. By the \Ito formula, we have that for any $t \in [0,T]$,
\begin{align}\label{eq17_1}
\E |e_\D (t \we \rho_R) |^2
& = \E \int_0^{t \we \rho_R}\Big ( 2 \lan X(s) - \yd (s), f(X(s),X(s-\d(s))) - \TDfds  \ran \no \\
& \qu + |g(X(s),X(s-\d(s))) - \TDgds|^2  \Big ) ds  \no \\
 & \le \E \int_0^{t \we \rho_R}\Big ( 2 \lan X(s) - \yd (s), f(X(s),X(s-\d(s))) - f(\yd (s), \yd (s- \d(s)))  \ran \no \\
& \qu + (p_1 - 1 )  |g(X(s),X(s-\d(s))) - g(\yd (s), \yd (s- \d(s)))|^2   \no \\
& \qu + 2 \lan X(s)-\yd (s), f(\yd (s), \yd (s- \d(s))) - \TDfds \ran \no \\
& \qu + \frac{p_1 -1}{p_1 -2}  |g(\yd (s), \yd (s- \d(s))) - \TDgds |^2   \Big ) ds  \no \\
&  \le  \Pi_1 + \Pi_2 ,
\end{align}
where
\begin{align*}
\Pi_1 : = &  \E \int_0^{t \we \rho_R} \Big ( |e_\D (s) |^2 + 2 \bar K_2 (|e_\D (s) |^2 + | X(s - \d(s)) - \yd (s - \d(s)) |^2 ) \no \\
& \qu - \frac{2}{1-\hat \d} U(X(s), \yd (s)) + 2U(X(s-\d(s)), \yd (s-\d(s))) \Big ) ds, \no \\
\Pi_2: = &  \E \int_0^{t \we \rho_R} \Big ( | f(\yd (s), \yd (s- \d(s))) - \TDfds |^2  \no \\
& \qu + \frac{p_1 -1}{p_1 -2}  |g(\yd (s), \yd (s- \d(s))) - \TDgds |^2   \Big ) ds ,
\end{align*}
where Assumption \ref{AS31_monotonity_and_Poly_cond} has been used. Noticing that $U(X(s), \yd (s)) = 0$ for any $s \in [-\tau,0 ]$, we then have
\begin{align*}
\int_0^{t \we \rho_R} U (X (s - \d(s)), \yd (s-\d(s))) ds  \le \frac{1}{1- \hat \d} \int_0^{t \we \rho_R} U (X (s), \yd (s)) ds
\end{align*}
and
\begin{align*}
\int_0^{t \we \rho_R} |X (s - \d(s))- \yd (s-\d(s))|^2  ds  \le \frac{1}{1- \hat \d} \int_0^{t \we \rho_R}|X (s )- \yd (s)|^2 ds.
\end{align*}
Consequently,
\begin{align}\label{eq17_2}
\Pi_1& \le \left (1+ 2 \bar K_1 +  \frac{2 \bar K_1}{1- \hat \d} \right ) \E \int_0^{t \we \rho_R} | e_\D (s) |^2 ds
\le \left (1+ 2 \bar K_1 +  \frac{2 \bar K_1}{1- \hat \d} \right ) \int_0^{t }\E | e_\D (s \we \rho_R) |^2 ds .
\end{align}
In order to estimate $\Pi_2$, we have the
\begin{align*}
\Pi_2  \le \Pi_{21} + \Pi_{22},
\end{align*}
where
\begin{align*}
\Pi_{21} & = 2 \int_0^T \E  | f(\yd (s), \yd (s-\d(s))) -  \fd(\yd (s), \yd (s-\d(s))) |^2 ds \no \\
& \qu + \frac{2(p_1 -1)}{p_1 -2} \int_0^T \E  | g(\yd (s), \gd (s-\d(s))) -  \fd(\yd (s), \yd (s-\d(s))) |^2 ds,
\end{align*}
and
\begin{align*}
\Pi_{22} & = 2 \int_0^T \E  | \fd (\yd (s), \yd (s-\d(s))) -  \TDfds |^2 ds \no \\
& \qu + \frac{2(p_1 -1)}{p_1 -2} \int_0^T \E  | \gd (\yd (s), \gd (s-\d(s))) -  \TDgds  |^2 ds.
\end{align*}
Noting that  $ p_0 \ge 2 + 3 l > 2 + l$,
in a similar way/fashion as Estimate (3.15) from \cite[Theorem 3.6]{Fei2020CPAA} was proved,
we also can show that
\begin{align}\label{eq18_1}
\Pi_{21} \le C [\mu^{-1} (\varphi (\Delta)) ]^{-(p_0 - l -2)}.
\end{align}
Recall \cite[Lemma 3.3]{Fei2020CPAA} that
 \begin{align*}
&|\fd(x_1,y_1) - \fd (x_2,y_2)|^2
 \le \bar K_2(1+|x_1|^l + |x_2|^l + |y_1|^l + |y_2|^l)   (|x_1-x_2|^2 + |y_1-y_2|^2),
\end{align*}
 for any $ x_1,y_1,x_2,y_2 \in \R^d$.
By the condition that $ p_0 \ge 2 + 3 l$ which implies $\displaystyle \frac{2p_0}{p_0 - l} \le \frac{p_0}{ 1 + l/2}$,
 \eqref{eq16_2} and the \Holder inequality as well as Lemma \ref{Lem3.4}, we have
 that for $s \in [0,T]$,
\begin{align}\label{eq18_3}
& \E | \fd (\yd (s), \yd (s-\d(s))) -  \TDfds |^2  \no \\
& \le \bar K_2 \E \Big ( \left [ |\yd (s) - Z_1(s)|^2 + | \yd (s- \d(s)) - Z_2(s)  |^2 \right ] \no \\
& \qu \K \left [ 1 + |\yd (s)|^l + |\yd (s - \d(s))|^l + |Z_1(s)|^l+ |Z_2(s)|^l \right ] \Big ) \no \\
& \le C \left ( \E | \yd (s) - Z_1(s)|^{2p_0/(p_0 - l)}+  \E | \yd (s-\d(s)) - Z_2(s)|^{2p_0/(p_0 - l)}    \right )^{(p_0 - l )/ p_0} \no \\
& \qu \K \left ( 1 + \E |\yd (s) |^{p_0} +  \E |\yd (s- \d(s) ) |^{p_0} + \E |Z_1 (s) |^{p_0} + \E |Z_2 (s) |^{p_0}\right )^{l/p_0} \no \\
& \le  C \left ( \E | \yd (s) - Z_1 (s)|^{2p_0 / (p_0 - l)}  \right )^{(p_0 - l )/ p_0}
+ C \left ( \E | \yd (s-\d(s)) - Z_2 (s)|^{2p_0 / (p_0 - l)}  \right )^{(p_0 - l )/ p_0} \no \\
&\le C  \D^{1 \we 2 \varrho }.
\end{align}
Similarly, we can show that
\begin{align}\label{eq20_1}
\E | \gd (\yd (s), \yd (s-\d(s))) -  \TDgds |^2 \le C  \D^{1 \we 2 \varrho }.
\end{align}
Combining \eqref{eq17_1}-\eqref{eq20_1}, we get
\begin{align}\label{eq20_3}
\E |e_\D(t \we \rho_R)|^2 \le C \intt  \E |e_\D(s \we \rho_R)|^2 ds + C \left ( \D^{2 \varrho \we 1 } \ve [\mu^{-1} (\varphi (\Delta))]^{-(p_0 - l -2)} \right ).
\end{align}
The Gronwall inequality give
\begin{align}\label{eq20_2}
\E |e_\D(t \we \rho_R)|^2  \le  C \left ( \D^{2 \varrho \we 1}  \ve [\mu^{-1} (\varphi (\Delta))]^{-(p_0 - l -2)} \right ).
\end{align}
Letting $R \to \infty$ gives assertion \eqref{eq17_5}.
\eqref{eq17_6} follows from \eqref{eq17_5} and Lemma \ref{Lem3.4}.
Finally, from \eqref{eq16_4}  and \eqref{eq21}, we may
define $\mu(R)$ and $\varphi(\D)$  by \eqref{phy11}, e.g.,
\begin{align*}
\mu (r) = \bar K_3 ^{1/2} r^{l/2}, \; \f r \ge 1 \qu \textrm{and}\qu \varphi(\D) = \hat h \D^{-1/4}, \; \f \D \in(0,1].
\end{align*}
Then
\begin{align}\label{temp104}
[ \mu^{-1} (\varphi (\Delta))  ]^{-(p_0 - l -2)} = C \D ^{ 0.5(p_0 -l- 2)/l} \le C \D ,
\end{align}
due to  $ p_0 \ge 2 + 3l$, which implies that
$0.5(p_0 -l- 2)/ l \ge 1$.
From \eqref{temp104} and \eqref{eq17_5} as well as \eqref{eq17_6}, we  obtain the  assertion \eqref{asser31}.
Thus, the proof is complete. $\Box$


%
\section{Mean-square  and $H_\infty$ stabilities }
In this section, we mainly discuss the mean-square  and $H_\infty$ stabilities of the truncated EM method for SDDE \eqref{eq0}. Noting that the truncated functions $f_\D$ and $g_\D$ can preserve the Khasminskii-type condition \eqref{PP3}, unfortunately  they cannot preserve the stability condition. Thus,
 we can only hope that the terms that work for the stability in the coefficients grows at most linearly,  while for those that grow super-linearly,
they have no stabilizing effect.  In our truncated method for stability, we only use the  truncation technique to
the super-linear terms in the coefficients. Lemma \ref{Lem22_conserving_Kcond_stable}  show that those partially truncated functions
 defined by \eqref{eqq11}
 have the property of preserving the stability condition.

As a result,
we assume that $f$ and $g$ can be decomposed as
$f(x,y) = F_1(x,y) + F(x,y)$ and $g(x,y) = G_1(x,y) + G(x,y)$, where
$F_1, F : \R^d \K \R^d  \to \R^d $ and $G_1, G: \R^d \K \R^d \to \R^{d\K m}$.
Moreover,  $F_1(0,0) = F(0,0) = G_1(0,0) = G(0,0)=0$, the coefficients $F_1$, $F$, $G_1$, $G$ satisfy the following conditions.
\begin{assumption}\label{AA1'_hybirc_Lips}
For any $R >0$, there exists  constants $\bar L$ and $\bar L_R$ depending on $R$   such that
\begin{align}\label{eeq1}
 & | F_1(x_1,y_1) - F_1(x_2 ,y_2 )|^2 \ve  | G_1(x_1,y_1) - G_1(x_2 ,y_2 )|^2 \le \bar L (|x_1 - x_2 |^2 + |y_1 - y_2 |^2),
 \end{align}
  for any $  x_1, x_2, y_1, y_2 \in \R^d$ and
 \begin{align}\label{eeq2}
 & | F(x_1,y_1) - F(x_2 ,y_2 )|^2 \ve  | G(x_1,y_1) - G(x_2 ,y_2 )|^2  \le \bar L_R (|x_1 - x_2 |^2 + |y_1 - y_2 |^2),
\end{align}
 for any $  x_1, x_2, y_1, y_2 \in \R^d$ with  $|x_1| \ve |x_2| \ve | y_2| \ve |y_2| \le R$.
\end{assumption}

\begin{assumption}\label{AA3'_stable_cond}
There exist nonnegative constants $\th$, $\l_1$, $\l_2$, $\a_1$, $\a_2$, $\a_3$, $\a_4$  and  $ \b >2$  such that
\begin{align} \label{eq7_2}
& 2\lan x, F_1(x,y) \ran + (1 + \th) |G_1(x,y)|^2  \le - \l_1|x|^2 + \l_2 |y|^2, \; \f x,y \in \R^d,  \no  \\
& 2 \langle x,   F(x,y)\rangle  + (1+\th^{-1})  |G (x,y)|^2  \le \a_1|x|^2 + \a_2|y|^2 - \a_3 |x|^\b + \a_4 |y|^\b,\; \f x,y \in \R^d.
\end{align}
\end{assumption}
When $\th=0$, we set $\th^{-1}|G(x,y)|^2 = 0$, when $\th = \infty$, we set $\th |G_1(x,y)|^2=0$. Clearly, Assumption \ref{AA3'_stable_cond} implies that
\begin{align}\label{eq7_3}
2 \lan x, f(x,y) \ran + |g(x,y)|^2 \le - (\l_1 - \a_1)|x|^2 + (\l_2 + \a_2)|y|^2 - \a_3|x|^\b + \a_4 |y|^\b,\; \f x,y \in \R^d.
\end{align}
We conclude from  \citet[Theorem 3.6]{Song2013DCDS} that the SDDE \eqref{eq0} is stable in mean square sense, which can
be stated by the following lemma.

\begin{lemma}\label{Lem21_stable_xx}
Let Assumptions  \ref{AA1'_hybirc_Lips}, \ref{AA3'_stable_cond} and  \ref{A2_bar_delta}  hold with
\begin{align}\label{eq8_0}
\l_1 > \a_1 + \frac{1}{4}\l_2+ ( \lfloor  (1- \hat \d)^{-1} \rfloor +1 ) (\l_2 + \a_2) \qu \textrm{and} \qu \a_3 > ( \lfloor  (1- \hat \d)^{-1} \rfloor +1 )  \a_4 \ge 0,
\end{align}
where $( \lfloor  (1- \hat \d)^{-1} \rfloor +1 ) = \lf 1/(1 - \hat \d) \rf + 1$.
Then for any given initial data \eqref{eq3}, the unique global solution $X(t)$ to \eqref{eq0} has the property that
\begin{align}\label{eq8_1}
\limsup_{t \to \infty} \frac{\log \E |X(t)|^2}{ t} \le - \left ( \g^{\star} \we \frac{1}{\tau} \log \frac{\a_3}{ ( \lfloor  (1- \hat \d)^{-1} \rfloor +1 ) \a_4}  \right ) ,
\end{align}
and
\begin{align}\label{eq8_2}
\int_0^\infty \E |X(t)|^2 dt < \infty,
\end{align}
where $\g^{\star} >0$ is the unique root to the following equation
\begin{align}\label{eq8_3}
\l_1  =  \Big ( \a_1 + \frac{1}{4}\l_2 \Big ) +  ( \lfloor  (1- \hat \d)^{-1} \rfloor +1 ) (\l_2 + \a_2) \textrm{e}^{\g^{\star} \tau}+\g^{\star}.
\end{align}
\end{lemma}
\begin{lemma}\label{Lem22_conserving_Kcond_stable}
Let Assumptions \ref{AA1'_hybirc_Lips}, \ref{AA3'_stable_cond} and  \ref{A2_bar_delta}  hold with
\begin{align}\label{eq10_1}
\l_1 > \a_1 + \frac{1}{4}\l_2+ ( \lfloor  (1- \hat \d)^{-1} \rfloor +1 ) (\l_2 + \a_2) \qu \textrm{and} \qu \a_3 > ( \lfloor  (1- \hat \d)^{-1} \rfloor +1 )  \a_4 \ge 0.
\end{align}
For a given step size $\D \in (0,1]$,
define
\begin{align}\label{eqq11}
\fd (x,y) = F_1 (x,y) + F_\D (x,y) \; \textrm{and} \; \gd (x,y) = G_1(x,y) + G_\D (x,y),\; \f x,y \in \R^d,
\end{align}
where $F_\D (x,y) = F(\pd (x), \pd (y)) $ and $G_\D (x,y) = G(\pd (x), \pd (y)) $.
Then
\begin{align}\label{eq10_2}
2 \lan x, f_\D(x,y) \ran + |g_\D(x,y)|^2 & \le - (\l_1 - \a_1 -\frac{1}{4}\a_2 )|x|^2 + (\l_2 + \a_2)|y|^2 \no \\
& \qu - \a_3|\pd (x) |^\b + \a_4 |\pd (y)|^\b, \; \f x,y \in \R^d,
\end{align}
 and
 \begin{align}\label{eq11_3}
|\fd(x, y)|^2 \D  \le \e_\D ( |x|^2 + |y|^2),\; \f x,y \in \R^d,
\end{align}
where
$\e_\D  = (4 \bar L + 2 \bar L_1) \D + 8  (\varphi(\D))^2 \D$.
\end{lemma}
\pr
For any $x\in \R^d$ with $|x| \le \mu^{-1}(\varphi(\D))$ and any $y \in \R^d$,
\eqref{eq10_2} follows from Assumption \ref{AA3'_stable_cond}. While for any $x\in \R^d$ with $|x| > \mu^{-1}(\varphi(\D))$ and any $y \in \R^d$, in a similar was as
\citet[Inequality (3.4)]{Fei2020CPAA} was obtained, we also derive  from Assumption \ref{AA3'_stable_cond} that
\begin{align}\label{eq9_4}
2 \lan x, f_\D(x,y) \ran + |g_\D(x,y)|^2 & \le (- \l_1 |x|^2 + \l_2 |y|^2)  - \a_3|\pd (x) |^\b + \a_4 |\pd (y)|^\b \no \\
& \qu + \frac{|x|}{ \mu ^{-1}(h (\D))} \left (  \a_1 |\pd (x)|^2+ \a_2 |\pd(y)|^2 \right )  \no \\
& \le (- \l_1 |x|^2 + \l_2 |y|^2)  - \a_3|\pd (x) |^\b + \a_4 |\pd (y)|^\b \no \\
& \qu + \a_1|x|^2 + \a_2|x||y|.
\end{align}
Now, using the inequality $\displaystyle |x||y| \le \frac{1}{4}|x|^2 +  |y|^2$ for any $x,y \in \R^d$ and rearranging the RHS terms
 of the last inequality in \eqref{eq9_4}, we obtain  \eqref{eq10_2}.
Now  let us establish  \eqref{eq11_3}
From Assumption \ref{AA1'_hybirc_Lips} and the property of truncated function, we have that
\begin{align}\label{eq12'_1}
& |\fd (x,y)|^2  = |F_1(x,y) + F_\D (x,y)|^2  \no \\
& \le 2 |F_1(x,y)|^2 + 2 |F_\D (x,y)|^2
\le 2 (\bar L + \bar L_1)(|x|^2 + |y|^2),  \; \textrm{if} \; |x|\le 1, \; |y| \le 1,
\end{align}
and
\begin{align}\label{eq12'_2}
|\fd (x,y)|^2
& \le 2 \bar L (|x|^2 + |y|^2) + 2(\varphi(\D))^2(1 + |x| + |y|)^2  \no \\
&\le ( 2 \bar L + 8 (\varphi(\D))^2 ) (|x|^2 + |y|^2),  \; \textrm{if} \; |x|\le 1, \; |y| > 1.
\end{align}
This also holds for $|x| >1$, $|y| \le 1$ or $|x| \ge 1$, $|y| \ge 1 $.
Thus, the proof is finished.
 $\Box$
 \par
The following theorem shows that the partially truncated EM solution can share the mean-square and $H_\infty$ stabilities of the true solution.
\begin{theorem}\label{Th_stable_Y}
Let Assumptions \ref{AA1'_hybirc_Lips}, \ref{AA3'_stable_cond} and  \ref{A2_bar_delta}   hold with
\begin{align*}
\l_1 > \a_1 + \frac{1}{4}\l_2+ ( \lfloor  (1- \hat \d)^{-1} \rfloor +1 ) (\l_2 + \a_2) \qu \textrm{and} \qu \a_3 > ( \lfloor  (1- \hat \d)^{-1} \rfloor +1 )  \a_4 \ge 0.
\end{align*}
Choose  $\D^{\star} \in (0,1]$ satisfying
\begin{align}\label{eg57}
\e_{\D^{\star} } =(4 \bar L + 2 \bar L_1) \D^\star + 8  (h(\D^\star))^2 \D^\star
= \frac{\l_1 - \a_1 - \frac{1}{4}\l_2 - ( \lfloor  (1- \hat \d)^{-1} \rfloor +1 ) (\a_2 +\l_2) }{{ 1 +  ( \lfloor  (1- \hat \d)^{-1} \rfloor +1 )  }}.
\end{align}
Then for any $\D \in (0, \D^{\star})$ and  any initial data \eqref{eq3}, the truncated EM approximation $ \{ y_k\}_{k \ge 0 }$ with
the truncated coefficients $\fd$ and $\gd$ given by \eqref{eqq11} has
the property that
\begin{align}\label{eq14_1}
 \limsup_{k \to \infty} \frac{\log \E |y_k|^2}{\tk} \le - \left ( \g^{\star }_\D  \we \frac{1}{\tau} \log \frac{\a_3}{( \lfloor  (1- \hat \d)^{-1} \rfloor +1 ) \a_4} \right ),
\end{align}
where $\g^{\star}_\D$ is the unique root to the following equation
\begin{align}\label{eq14_2}
\l_1 = \Big  ( \a_1 + \frac{1}{4}\l_2 + \e_\D  \Big ) + ( \lfloor  (1- \hat \d)^{-1} \rfloor +1 ) (\l_2 + \a_2 + \e_\D)
\textrm{e}^{ \g_\D ^{\star}\tau} + \frac{1 - \textrm{e}^{-\g_\D ^{\star} \D }  }{\D}.
\end{align}
Moreover,
\begin{align}\label{eq14_3}
\lim_{\D \to 0} \g^{\star}_\D= \g^{\star},
\end{align}
and
\begin{align}\label{eq133_6}
\int_0^{\infty} \E |Z_ 1(s)|^2 ds < \infty.
\end{align}
\end{theorem}
\pr
From  \eqref{Discrete_EM}, we have
\begin{align}\label{eq11_1}
|y_{k+1}|^2 = |y_k|^2+ 2 \lan y_k, \fd(y_k, y_{k-\d_k}) \ran \D +
|\gd(y_k, y_{k-\d_k})|^2 \D  + | \fd(y_k, y_{k-\d_k})|^2 \D^2 + J_k,\; \forall k \ge 0,
\end{align}
where
\begin{align*}
J_k = 2 \lan y_k, \gd(y_k, y_{k-\d_k}) \D B_k \ran + 2 \lan \fd(y_k, y_{k-\d_k}), \gd(y_k, y_{k-\d_k}) \D B_k \ran\D
 + |\gd(y_k, y_{k-\d_k})|^2(|\D B_k|^2 -\D ).
\end{align*}
Obviously, $\E J_k =0$.
Using  Lemma \ref{Lem22_conserving_Kcond_stable} yields
\begin{align}\label{eq11_2}
& 2 \lan y_k, \fd(y_k, y_{k-\d_k}) \ran  + |\gd(y_k, y_{k-\d_k})|^2   \no \\
& \le - \Big (\l_1 - \a_1 - \frac{1}{4}\a_2 \Big )|y_k|^2 + (\l_2 + \a_2)|\ykd|^2
-\a_3|\pd (y_k)|^\b + \a_4 |\pd (\ykd)|^\b.
\end{align}
Inserting \eqref{eq11_2}  into \eqref{eq11_1} and using  \eqref{eq11_3}, we have
\begin{align}\label{eq11_4}
|y_{k+1}|^2 & \le |y_k|^2 - \Big (  \l_1 -\a_1 - \frac{1}{4} \a_2 - \e_\D  \Big )|y_k|^2 \D + ( \l_2 + \a_2 + \e_\D )  |\ykd|^2\D \no \\
& \qu -\a_3|\pd (y_k)|^\b \D + \a_4 |\pd (\ykd)|^\b \D + J_k,\; \forall k \ge 0,.
\end{align}
For an arbitrary constant $r >1$, we see from \eqref{eq11_4} that
\begin{align}\label{eq12_1}
& r^{(k+1)\D}\E|\ykk|^2  - r^{k\D} \E|\yk|^2 \no \\
& \le (r^{(k+1)\D} - r^{k \D })\E |\yk|^2 - \Big (  \l_1 -\a_1 - \frac{1}{4} \a_2 - \e_\D  \Big )  r^{(k+1)\D} \E |\yk|^2 \D
+( \l_2 + \a_2 + \e_\D )  r^{(k+1)\D} \E |\ykd|^2 \D \no \\
& \qu + \E \Big [  - \a_3 r^{(k+1)\D}  |\pd (\yk)|^\b\D + \a_4 r^{(k+1)\D} |\pd (y_{k-\d_k})|^\b \Big ]\D ,\; \forall k \ge 0, .
\end{align}
Consequently,
\begin{align}\label{eq12_2}
r^{(k+1)\D} \E |\ykk|^2 & \le |\x (0)|^2 + \left (- \Big (  \l_1 -\a_1 - \frac{1}{4} \a_2 - \e_\D  \Big ) \D + 1 - r^{-\D}  \right ) \sum_{j=0}^k r^{(j+1)\D} \E |y_j|^2  \no \\
& \qu  +( \l_2 + \a_2 + \e_\D )  \sum_{j=0}^k r^{(j+1)\D} \E |y_{j-\d_j}|^2  \D  \no \\
& \qu + \E \left  [ - \a_3   \sum_{j=0}^k r^{(j+1)\D}  | \pd (y_j)|^\b
   + \a_4  \sum_{j=0}^k r^{(j+1)\D}  | \pd (y_{j-\d_j})|^\b \right  ] \D,\; \forall k \ge 0, .
\end{align}
By Lemma \ref{Lem_M}, we get
\begin{align}\label{eq11_5}
\sum_{j=0}^k r^{(j+1)\D}  |y_{j-\d_j}|^2 & \le  ( \lfloor  (1- \hat \d)^{-1} \rfloor +1 ) r^{\tau}  \sum_{j=-M}^k r^{(j+1)\D}  |y_{j}|^2 \no \\
& = ( \lfloor  (1- \hat \d)^{-1} \rfloor +1 ) r^{\tau}  \sum_{j=-M}^{-1} r^{(j+1)\D}  |y_{j}|^2 + ( \lfloor  (1- \hat \d)^{-1} \rfloor +1 ) r^{\tau}  \sum_{j=0}^{k} r^{(j+1)\D}  |y_{j}|^2 \no \\
& \le \frac{( \lfloor  (1- \hat \d)^{-1} \rfloor +1 ) r^{\tau}}{ 1 - r^{-\D}} \| \x \|^2 +  \bar  \k r^{\tau} \sum_{j=0}^k r^{(j+1)\D}  |y_{j}|^2,\; \forall k \ge 0,
\end{align}
and
\begin{align}\label{eq12_3}
\sum_{j=0}^k r^{(j+1)\D}  |\pd (y_{j-\d_j})|^\b
& \le \frac{( \lfloor  (1- \hat \d)^{-1} \rfloor +1 ) r^{\tau}}{ 1 - r^{-\D}} \| \x \|^\b +  \bar  \k r^{\tau}  \sum_{j=0}^k r^{(j+1)\D}  |\pd(y_{j})|^\b,\; \forall k \ge 0.
\end{align}
Inserting \eqref{eq11_5} and \eqref{eq12_3} into \eqref{eq12_2} gives that
\begin{align}\label{eq12_5}
r^{(k+1)\D} \E |\ykk|^2 & \le H_0 (r,\D) - H_1 (r,\D)  \sum_{j=0}^k r^{(j+1)\D} \E |y_j|^2 \D \no \\
& \qu - H_2 (r) \sum_{j=0}^k r^{(j+1)\D} \E |\pd (y_j)|^2  \D , \; \forall k \ge 0,
\end{align}
where
\begin{align}\label{H_(Delta)}
H_0(r,\D) & = \| \x \|^2 + ( \lfloor  (1- \hat \d)^{-1} \rfloor +1 ) r^\tau \Big [ (\l_2 + \a_2 + \e_\D) \|  \x \|^2 +  \a_4 \|  \x \|^\b \Big ] \frac{\D}{ 1 - r^{-\D}}, \no \\
H_1(r,\D) & = \Big [\l_1 -\a_1 - \frac{1}{4} \a_2 - \e_\D - ( \lfloor  (1- \hat \d)^{-1} \rfloor +1 ) r^{\tau} (\l_2 + \a_2 + \e_\D)\Big ] - \frac{1- r^{-\D}}{\D}, \no \\
H_2(r) & = \a_3 - ( \lfloor  (1- \hat \d)^{-1} \rfloor +1 ) r^{\tau } \a_4.
\end{align}
Choose  $\D^{\star} \in (0,1]$  such that \eqref{eg57} holds, i.e.,
\begin{align*}
\e_{\D^{\star}} = \frac{\l_1 - \a_1 - \frac{1}{4}\l_2 - ( \lfloor  (1- \hat \d)^{-1} \rfloor +1 ) (\a_2 +\l_2) }{ 1 +  ( \lfloor  (1- \hat \d)^{-1} \rfloor +1 ) }.
\end{align*}
Then for any $\D \in (0, \D^{\star})$, we have
\begin{align}\label{eq13_4}
H_1(1,\D) =  \l_1 - \a_1 - \frac{1}{4}\l_2 - \e_\D - ( \lfloor  (1- \hat \d)^{-1} \rfloor +1 ) (\l_2 + \a_2 + \e_\D) >0
\end{align}
\begin{align}\label{eq13_5}
H_1(\bar r, \D) = - \frac{1 - \bar r ^{-\D}}{ \D} <0, \; \textrm{with}\;  \left ( \frac{\l_1 - \a_1 - \frac{1}{4}\l_2 - \e_\D }{( \lfloor  (1- \hat \d)^{-1} \rfloor +1 ) (\l_2 + \a_2 + \e_\D)} \right )^{1/\tau}  >1
\end{align}
and \begin{align}\label{eq13_6}
\frac{d H_1 (r,\D)}{ dr} <0.
    \end{align}
From \eqref{eq13_4}, \eqref{eq13_5} and \eqref{eq13_6},  there is a positive constant $r_1^{\star} =r_1^{\star}(\D) \in (1, \bar r)$
such that $H_1(r_1^{\star}, \D) =0$. Let $r^{\star} = r^{\star}(\D) = r_1^{\star}(\D)\we r_2 ^{\star}$ with
$ \displaystyle
r_2^{\star} =
\left ( \frac{\a_3}{ ( \lfloor  (1- \hat \d)^{-1} \rfloor +1 ) \a_4}  \right )^{1/\tau} >1
$,
then for any $1 < r < r^{\star}(\D)$, we have $H_1(r,\D) >0$ and $H_2(r)>0$, thus we conclude from \eqref{eq12_5} that
\begin{align*}
r^{(k+1)\D} \E |\ykk|^2 \le H_0 ( r,\D)< \infty, \;\f k \ge 0.
\end{align*}
Therefore,
\begin{align}\label{eq13_8}
\limsup_{k \to \infty} \frac{\log \E |\yk|^2}{ \tk} \le - \log r.
\end{align}
Bearing in mind that $r_2^{\star} = \textrm{e}^{\g^{\star}_2 } $ and
setting $r = \textrm{e}^{\g}$, $r_1^{\star} = \textrm{e}^{\g_1 } $, then
\eqref{eq13_8} becomes \eqref{eq14_1}.
Finally, noticing that $\e_\D \to 0$ and $\displaystyle \frac{1 - \textrm{e}^{- r_1 ^{\star} \D}}{\D } \to r_1^{\star} $
as $\D \to 0$, comparing \eqref{eq8_3} with \eqref{eq14_2}, we obtain the assertion \eqref{eq14_3}.  \\
Finally, we begin to establish \eqref{eq133_6}.
By \eqref{eq11_4}, we have
\begin{align}\label{eq13'_1}
\E |\ykk|^2 & \le |\x (0)|^2 -\Big (  \l_1 -\a_1 - \frac{1}{4} \a_2 - \e_\D  \Big )  \sum_{j=0}^k  \E|y_j|^2 \D+ ( \l_2 + \a_2 + \e_\D )  \sum_{j=0}^k \E |y_{j-\d_j}|^2  \D \no \\
& \qu + \E \left [ - \a_3  \sum_{j=0}^k  |\pd (y_j)|^\b  +\a_4  \sum_{j=0}^k  |\pd (y_{j- \d_j})|^\b \right ] \D, \; \forall k \ge 0,.
\end{align}
While Lemma \ref{Lem_M} gives that
\begin{align}\label{eq13'_2}
\sum_{j=0}^k |y_{j-\d_j}|^2 \le ( \lfloor  (1- \hat \d)^{-1} \rfloor +1 ) M \| \x \|^2+ \sum_{j=0}^k |y_{j}|^2
\end{align}
and
\begin{align}\label{eq13'_3}
\sum_{j=0}^k |\pd ( y_{j-\d_j})|^\b  \le ( \lfloor  (1- \hat \d)^{-1} \rfloor +1 ) M \| \x \|^\b+ \sum_{j=0}^k |\pd (y_{j})|^\b.
\end{align}
Substituting \eqref{eq13'_2} and \eqref{eq13'_3} into \eqref{eq13'_1}, we obtain that for any $\D \in (0, \D^{\star})$,
\begin{align}
0 \le \E|\ykk|^2 & \le \| \x \|^2 + ( \lfloor  (1- \hat \d)^{-1} \rfloor +1 ) \tau \Big ( ( \l_2 + \a_2 + \e_\D ) \| \x \|^2 + \a_4 \| \x \|^\b \Big ) \no \\
& \qu - (\a_3 - ( \lfloor  (1- \hat \d)^{-1} \rfloor +1 ) \a_4) \E\left [  \sum_{j=0}^k|\pd (y_{j})|^\b \right ] \D \no \\
& \qu - \Big [ \l_1 -\a_1 - \frac{1}{4} \a_2 - \e_\D - ( \lfloor  (1- \hat \d)^{-1} \rfloor +1 )  (\l_2 + \a_2 + \e_\D)\Big ]
   \sum_{j=0}^k  \E| y_j |^2 \D  , \; \forall k \ge 0,
\end{align}
which means
\begin{align}
\sum_{j=0}^k \E|\pd (y_{j})|^2 \D  \le \frac{\| \x \|^2 + ( \lfloor  (1- \hat \d)^{-1} \rfloor +1 ) \tau \Big (  ( \l_2 + \a_2 + \e_\D ) \| \x \|^2 + \a_4 \| \x \|^\b \Big )}{\l_1 -\a_1 - \frac{1}{4} \a_2 - \e_\D - ( \lfloor  (1- \hat \d)^{-1} \rfloor +1 )  (\l_2 + \a_2 + \e_\D) } < \infty,
\end{align}
holds for  any $k \ge 0$.
Letting $k \to \infty$ gives \eqref{eq133_6}.
Thus, the proof is finished.
$\Box$

\section{Numerical examples}

\begin{example}\label{example1}
Consider a highly nonlinear scalar SDDE (see \cite{Fei2020CPAA})
\begin{align}\label{eg1}
dX(t)&  = \left [ -9X^3 (t) + |X (t - \d(t))|^{3/2} \right  ]dt + X^2(t) d B(t), \; t \ge 0,
\end{align}
with initial data $\{ X(t): - \tau \le t \le 0 \}= \xi  \in C([-\tau,0]; \R)$,
where $B(t)$ is a one-dimensional Brownian motion.
Assume that
$\d $ satisfies Assumption \ref{A2_bar_delta}.
Clearly, the coefficients
\begin{align}\label{eeg7}
f(x,y) = -9 x^3 + |y|^{3/2} \qu \textrm{and} \qu g(x,y) = x^2 , \; \f x,y \in \R,
\end{align}
are locally Lipschitz continuous. Moreover, if $p_0 = 18.5$, then
\begin{align}\label{eeg1}
& x f(x,y) + \frac{p_0 - 1}{2 }  = -9x^4 + x|y|^{3/2} + 8.75x^4 \no \\
& \le -9x^4  + 8.75x^4 + 0.25x^4 + 0.75y^2 = 0.75y^2,
\end{align}
which means that Assumption \ref{AA3_Kha_cond_P0} is satisfied. For any
$x_1, x_2, y_1, y_2, \in \R$, we have
\begin{align}\label{eeg2}
& (x_1 - x_2) (f(x_1, y_1) - f(x_2,y_2)) \no \\
& \le -4.5 (x_1^2 + x_2^2 )|x_1 - x_2|^2 + 0.5 |x_1 - x_2|^2
+ 2.25 |y_1 - y_2 |^2 + 2.25(y_1^2 + y_2^2)|y_1 - y_2|^2
\end{align}
 and
\begin{align}\label{eeg3}
|g(x_1,y_1) - g(x_2,y_2)|^2 = | x_1^2 - x_2^2|^2 \le 2 (x_1^2 + x_2^2 )|x_1 - x_2|^2,
\end{align}
see \cite[p. 2086]{Fei2020CPAA}.
Consequently,
\begin{align}\label{eeg4}
& (x_1 - x_2) (f(x_1, y_1) - f(x_2,y_2)) + \frac{p_1 - 1}{2} |g(x_1, y_1) - g(x_2,y_2)|^2 \no \\
& \le   0.5 |x_1 - x_2|^2 + 2.25 |y_1 - y_2 |^2
-(5.5- p_1) (x_1^2 + x_2^2 )|x_1 - x_2|^2
+ 1.125(y_1^2 + y_2^2)|y_1 - y_2|^2.
\end{align}
If we set $p_1 = 3.25$, $\hat \d = 0.5$ and $U(x_1, x_2 ) = 1.125(x_1^2 + x_2^2)|x_1 - x_2|^2$,
then \eqref{eeg4} becomes
\begin{align}\label{eeg5}
& (x_1 - x_2) (f(x_1, y_1) - f(x_2,y_2)) + \frac{p_1 - 1}{2} |g(x_1, y_1) - g(x_2,y_2)|^2 \no \\
& \le   0.5 |x_1 - x_2|^2 + 2.25 |y_1 - y_2 |^2
- \frac{1}{ 1 - \hat \d} U(x_1, x_2 ) + U(y_1, y_2 ).
\end{align}
Moreover,
it is straightforward to show that \eqref{eq16_2} is satisfied with $l =4$. Thus,
we have verified  Assumption \eqref{AS31_monotonity_and_Poly_cond} with $p_0  \ge 2  + 3 l$. From \eqref{eeg7} and \eqref{eq21}, we may set
\begin{align}\label{eeq14}
\mu (R) = 10 R^2,  \; \f R \ge 1 \qu  \textrm{and} \qu \varphi(\D) =10 \D^{- 1/4}, \; \f \D \in (0,1].
\end{align}
Then, by Theorem \ref{Th_3__convergence_rate}, for any $\D \in (0,1]$ the truncated EM solution $\yd (t)$  will converge to the exact solution $X(t)$   in the sense that
\begin{align}\label{eeg8}
 \E |X(T) - \yd (T)|^2  \le C \D^{  }\qu \textrm{and} \qu  \E |X(T) - Z_1(T)|^2  \le C \D^{  }, \; \f T>0.
\end{align}
However, if constant delay is considered in SDDE \eqref{eg1}, we may also apply the truncated EM method from \cite{Fei2020CPAA} to \eqref{eg1} by setting
\begin{align}\label{eeq11}
\mu (R) = 10 R^3,  \; \f R \ge 1 \qu  \textrm{and} \qu \varphi(\D) = 10 \D^{- 1/5},\; \f \D \in (0,1],
\end{align}
due to
\begin{align}\label{eeq12}
\sup_{ |x| \ve |y| \le R} \Big ( |f(x,y)| \ve |g(x,y)| \Big ) \le 10 R^3, \; \f R \ge 1.
\end{align}
Thus according to  \cite[Corollary 3.7]{Fei2020CPAA},   for any $\D \in (0,1]$ the truncated EM solution will approximate the exact solution
in the sense that
\begin{align}\label{eeg13}
 \E |X(T) - \yd (T)|^2  \le C \D^{3/5}\qu \textrm{and} \qu  \E |X(T) - Z_1(T)|^2  \le C \D^{3/5}, \; \f T>0,
\end{align}
see \cite[p. 2086]{Fei2020CPAA}.
Comparing \eqref{eeg8} and \eqref{eeg13}, we conclude that our
scheme
establishes a better rate of convergence result than that of \citet{Fei2020CPAA} under almost the same conditions. Set $ \tau = 1$, $ \d(t) = 0.5 - 0.5 \sin(t)$ and $X(t) = 2$ for any $t \in [- \tau , 0]$. Truncated EM solution with step size $\D = 2^{-14}$ is taken as the replacement of the true solution.
The root of mean-square errors  with different step sizes $2^{-7}$, $2^{-8}, \cdots$  $2^{-11}$  for $500$ at time $T = 10$ simulations is illustrated in Fig. 1(a). A least square fit of errors yields the strong convergence order $0.5134$ and thus is close to the theoretical value $0.5$.
\end{example}

\begin{example}\label{example2}
Let us consider the stochastic delay power logistic model
(see e.g., \cite{Mao2004Delay_pop,Mao2004delay_LV})
\begin{align}\label{eg2}
dX(t) = X(t)[ a + b X (t - \d (t)) - X^2 (t) ] dt + c X(t)X(t - \d (t)) d B(t), \; t \ge 0
\end{align}
with initial data $\{ X(t): - \tau \le t \le 0 \}= \xi  \in C([-\tau,0]; \R)$,
where $B(t)$ is a one-dimensional Brownian motion and $a, b,c$ are all constants.
Assume that $\xi$ satisfies Assumption \ref{AA4_initial_cond} with $K_4 = 2$,  $\varrho = 0.5$ and
$\d $ satisfies Assumption \ref{A2_bar_delta}.
Set
\begin{align*}
 f(x,y) = F_1 (x,y) + F (x,y)\qu  \textrm{ and} \qu g (x,y) = G_1 (x,y) + G(x,y)
\end{align*}
 where
 \begin{align*}
 F_1 (x,y) = a x, \qu G_1(x,y) = 0,\qu  F(x,y) = bxy - x^2, \qu G (x,y) = c xy,
 \end{align*}
for any $x,y \in \R$.
Clearly, Assumption \ref{AA1'_hybirc_Lips} hold. Put $\theta = \infty$. Then
\begin{align*}
2 x F_1(x,y) + (1 + \th )|G_1 (x,y)|^2 = 2a x^2,
\end{align*}
and the elementary inequality yields that
\begin{align*}
& 2 x F(x,y) + (1 + \th^{-1} )|G (x,y)|^2  = 2x ( bxy - x^3) + (c xy)^2 \no \\
& \le 0.5 x^4 + 0.5 (2 by)^2 - 2 x^4 + 0.5x^4 + 0.5(c^2 y^2)^2
 =2 b^2 y^2 - x^4 + 0.5c^4 y^4.
\end{align*}
Thus,
\begin{align*}
 \lambda _ 1 = -2a, \qu  \l_2 = 0, \qu \a_1 =0, \qu \a_2 = 2b^2, \qu \a_3 = 1,\qu
\a_4 = 0.5c^4, \qu \beta = 4.
\end{align*}
%
If  we let
\begin{align}\label{eg43}
0.5c^4 ( \lfloor  (1- \hat \d)^{-1} \rfloor +1 ) <1 \qu \textrm{and} \qu -2a >  2  b^2  ( \lfloor  (1- \hat \d)^{-1} \rfloor +1 ),
\end{align}
which means that condition \eqref{eq10_1} is satisfied,
then by Lemma \ref{Lem21_stable_xx},
the exact solution $X(t)$ to SDDE \eqref{eg2} has the property that
\begin{align}\label{eg46}
\limsup_{t \to \infty} \frac{\log \E |X(t)|^2}{ t} \le - \left ( \g^{\star} \we \frac{1}{\tau} \log \frac{1}{  0.5c^4( \lfloor  (1- \hat \d)^{-1} \rfloor +1 )}  \right ),
\end{align}
 and $\displaystyle \int_0^\infty \E |X(t)|^2 dt < \infty,$
where $\g^{\star} >0$ is the unique root to the following equation
\begin{align}\label{eg45}
-2a =  \g^{\star} + 2 b^2 ( \lfloor  (1- \hat \d)^{-1} \rfloor +1 )   \textrm{e}^{\g^{\star} \tau}.
\end{align}
On the other hand, take
 \begin{align}\label{eeq69}
\mu (R) =  \Big ( ( |b| + 1)  \ve |c| \Big ) R^2,  \; \f R \ge 1 \qu  \textrm{and} \qu \varphi(\D) = \Big ( ( |b| + 1)  \ve |c| \Big ) \D^{- 1/4}, \; \f \D \in (0,1].
\end{align}
 We apply partially  truncated EM scheme \eqref{Discrete_EM} with coefficients $\fd$ and $\gd$ given by \eqref{eqq11} and we denote
$ \{ y_k\}_{k \ge 0 }$ by the discrete  truncated EM approximation. According to Theorem \ref{Th_stable_Y}, for any $\D \in (0, \D^{\star})$ and  any initial data, the
truncated EM approximation $ \{ y_k\}_{k \ge 0 }$ has the property that
\begin{align}\label{eg60}
 \limsup_{k \to \infty} \frac{\log \E |y_k|^2}{\tk}
 \le - \left ( \g^{\star }_\D  \we \frac{1}{\tau} \log \frac{1}{0.5 c^4 ( \lfloor  (1- \hat \d)^{-1} \rfloor +1 ) } \right ),
\end{align}
and $\displaystyle \sum_{k = 0}^{\infty} \E |y_k|^2 \D  < \infty,$
where $\g^{\star}_\D >0$ is the unique root to the following equation
\begin{align}\label{eg48}
-2a = \frac{1 - \textrm{e}^{-\g_\D ^{\star} \D }  }{\D} + \e_\D + ( \lfloor  (1- \hat \d)^{-1} \rfloor +1 ) ( 2 b^2 + \e_\D)     \textrm{e}^{\g^{\star} \tau}.
\end{align}
Numerically, if we set
\begin{align*}
\hat \d = 0.05, \qu \tau = 0.1, \qu  a = -3, \qu b = 1, \qu c = 0.5, \qu \varphi(\D) = 2 \D ^{-1/4},
\end{align*}
 then
 \begin{align*}
( \lfloor  (1- \hat \d)^{-1} \rfloor +1 ) = 2, \qu \g^{\star} = 1.3992, \qu \bar L = 5,  \qu \bar L_R = (3 R^2 + 1), \qu  \e_\D = 20 \D + 32 \D^{1/2},
 \end{align*}
and solving \eqref{eg57} gives $ \D ^{\star} = 4.2308 \K 10^{-4}$.
Computational results for
$\e_\D$ and $\g^{\star}_\D $ with different step sizes $\D$ are shown in Table \ref{table1}.
Fig illustrates a simple path of the truncated EM solution $\yd (t)$ with with step size $\D = 10^{-4}$ and
$\d (t) = 0.05- 0.05 \sin(t)$.

\end{example}

\begin{table}[!t]
 \caption{ $\epsilon_\D$ and $\g_\D^{\star}$ with different step sizes for solving \eqref{eg48}
 }
 \label{table1}
 \centering
 \begin{tabular}{ccccccc}
  \toprule
  $\D$  & $10^{-4}$ & $10^{-5}$  & $10^{-6}$  & $10^{-7}$ &  $10^{-8}$   &  $10^{-9}$    \\
  \midrule
$    \epsilon_\D   $    & 0.3220   & 0.1014 & 0.0320 &  0.0101   & 0.0032   & 0.0010  \\
  $\g_\D^{\star}  $    & 0.6982   & 1.1728 & 1.3272 &  1.3764   & 1.3920   &  1.3970 \\
  \bottomrule
 \end{tabular} \\
\end{table}
\begin{figure}[htb]
  \centering
  \label{fig1}
  \subfloat[Root of mean-square errors    for  \eqref{example1}]{
    \includegraphics[width=0.50\textwidth]{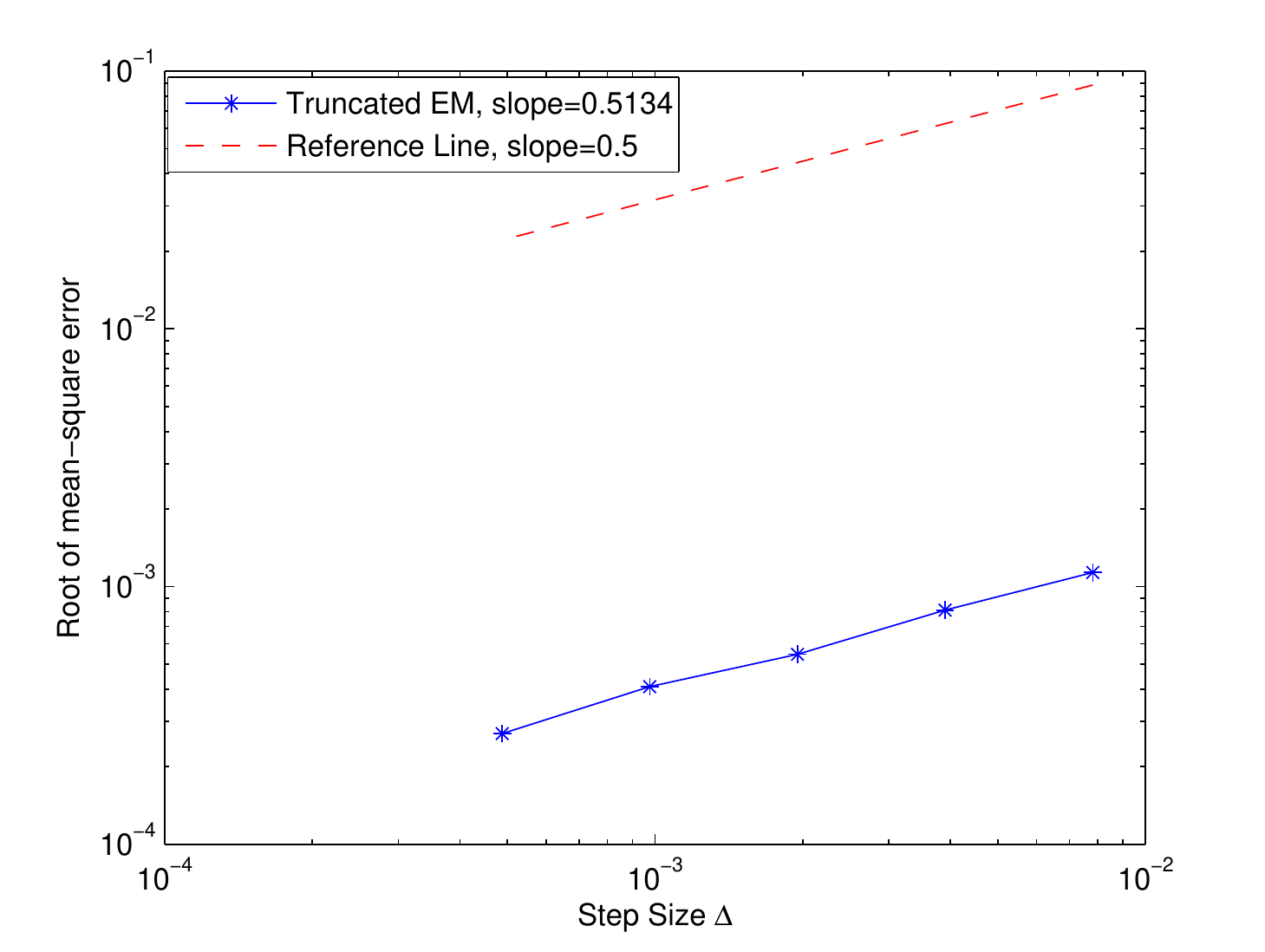}
  }
  \subfloat[A sample path of  $Y_\D (t)$    for  \eqref{eg2}]{
    \includegraphics[width=0.50\textwidth]{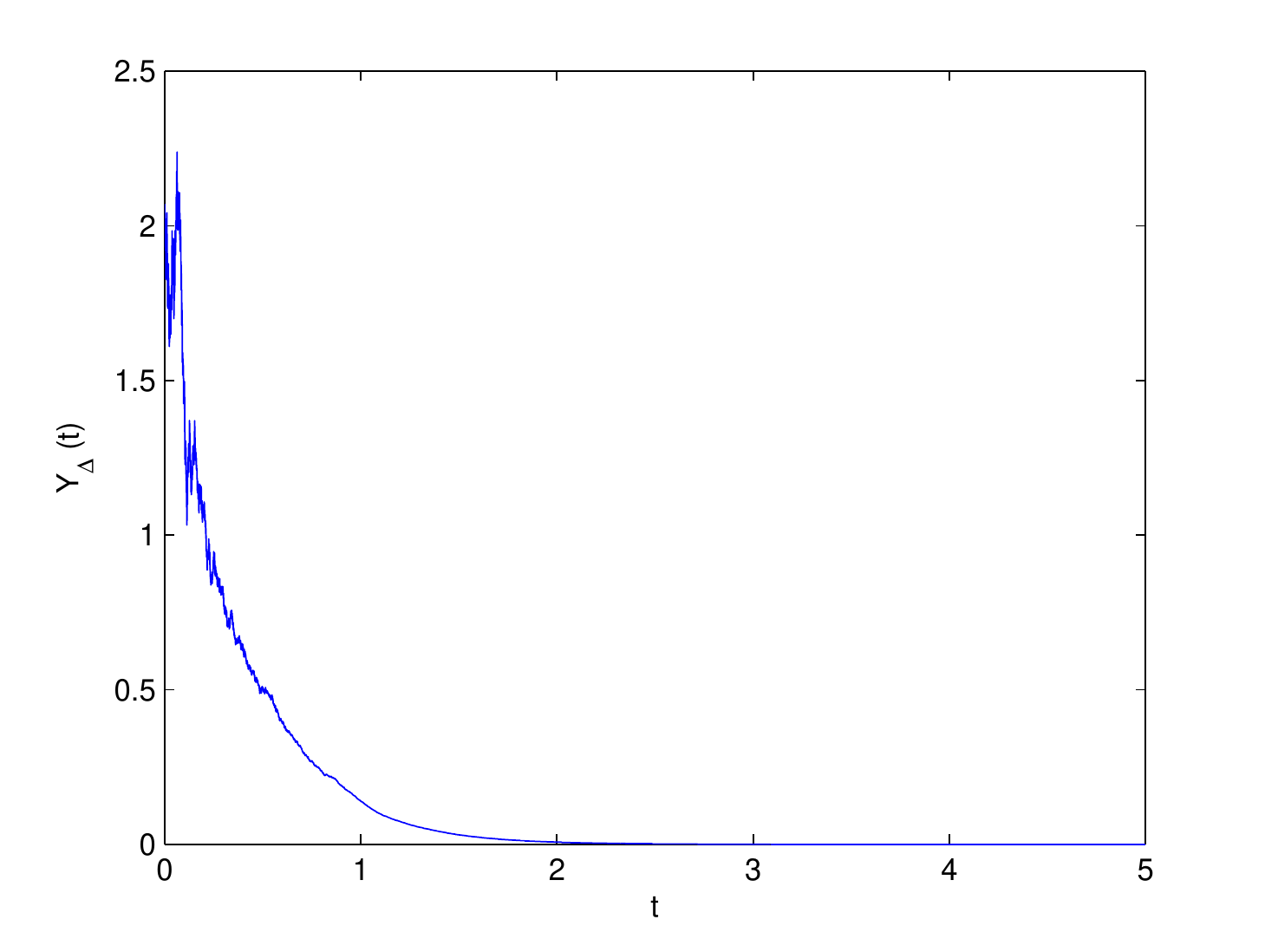}
  }
  \hfil
    \caption{Numerical simulations for \eqref{eg1} and  \eqref{eg2}}
\end{figure}

\section{Conclusion}

In this paper, we mainly study the strong convergence and stability of truncated EM scheme for SDDEs with variable delay.
The results show that our method has  a  better convergence order than those of the existing literature on the truncated EM for SDDEs under more relaxing conditions. Numerical simulations are provided to show the effectiveness of the  theoretical results.
In the future, we will consider  the strong convergence of the numerical scheme  for SDDEs driven
by L\'{e}vy process where all the three coefficients can grow super-linearly  .

\section*{Acknowledgments}
This work was supported in part by the Natural Science Foundation of China (No. 71571001).

\section*{References}

\end{document}